\documentclass[a4paper,leqno]{amsart}

\usepackage{amsmath}
\usepackage{amsthm}
\usepackage{amssymb}
\usepackage{amsfonts}
\usepackage[applemac]{inputenc}
\usepackage{mathrsfs}
\usepackage{pdfsync}
\usepackage{graphicx}
\usepackage{mathabx}
\usepackage{color}

\def\C{{\mathbb C}}
\def\R{{\mathbb R}}
\def\N{{\mathbb N}}
\def\Z{{\mathbb Z}}

\def\le{\leqslant}
\def\ge{\geqslant}
\def\eps{\varepsilon}
\def\Eq#1#2{\mathop{\sim}\limits_{#1\rightarrow#2}}

\theoremstyle{plain}
\newtheorem{theorem}{Theorem}[section]

\theoremstyle{definition}

\newtheorem{remark}[theorem]{Remark}
\newtheorem*{remark*}{Remark}

\numberwithin{equation}{section}

\begin{document}

\title[Singularity formation in 1D defocusing fNLS]
{Numerical evidence for singularity formation in defocusing fractional NLS in one space dimension}

\author[C. Klein]{Christian Klein}
\address[C.~Klein]
{Institut de Math\'ematiques de Bourgogne, Universit\'e de 
Bourgogne-Europe, Institut Universitaire de 
France,
9 avenue Alain Savary, BP 47870, 21078 Dijon Cedex}
\email{christian.klein@u-bourgogne.fr}

\author[C. Sparber]{Christof Sparber}
\address[C.~Sparber]
{Department of Mathematics, Statistics, and Computer Science, M/C 249, University of Illinois at Chicago, 851 S. Morgan Street, Chicago, IL 60607, USA}
\email{sparber@uic.edu}

\begin{abstract}
We consider nonlinear dispersive equations of Schr\"odinger-type involving fractional powers $0<s\le 1$ of the Laplacian and a defocusing power-law nonlinearity. 
We conduct numerical simulations in the case of small, energy supercritical $s$ and provide 
evidence for a novel type of highly oscillatory singularity within the solution. 
\end{abstract}

\date{\today}


\thanks{This publication is supported by the MPS Simons foundation 
through award no. 851720, the ANR project ANR-17-EURE-0002 EIPHI, and 
the ANR project ISAAC-ANR-23-CE40-0015-01.
}
\maketitle


\section{Introduction}
\label{sec:intro}

In this work, we shall be interested in fractional nonlinear Schr\"odinger equations (fNLS) of the form
\begin{equation}\label{fNLS}
i \partial_t u = (- \Delta )^s u + |u|^{2p} u, \quad u_{\vert{t=0}} = u_0,
\end{equation}
where $(t,x) \in \R\times \R^d$ (although we shall mainly study the case $d=1$), $p \in \N$ a nonlinear parameter, and $0<s\le 1$, the 
{\it fractional power} of the Laplacian. The latter is defined as a Fourier-multiplier with symbol $\lambda(\xi)=|\xi|^{2s}$, i.e. 
\[
(- \Delta )^s u = \mathcal F^{-1} \big(|\xi|^{2s} \widehat u \big),
\]
where $\widehat u\equiv \mathcal F u$ denotes the Fourier transform of $u$. 

Linear equations of the form \eqref{fNLS} were introduced by Laskin in an attempt to generalize the Feynman path integral of classical quantum mechanics to Levy processes \cite{La}. 
Their nonlinear counterpart appears in various applications, for example 
in the description of water waves \cite{IoPu}, or in the continuum limit of certain long-range lattice-models \cite{KLS}. 
Due to their strong resemblances to the well-studied classical NLS (obtained for $s=1$), they 
are an active field of mathematical research in which one tries to extend previous results 
to a fractional, and thus nonlocal setting, see, e.g., \cite{BHL, CA, HoSi, KMS, KLR} and the references therein.

Let us recall that solutions to \eqref{fNLS} formally satisfy the conservation laws of {\it mass}
\begin{equation}\label{mass}
M(u) = \int_{\R^d} |u|^2 \, dx,
\end{equation}
and {\it energy}
\begin{equation}\label{energy}
E(u) = \frac{1}{2} \int_{\R^d} \big| (- \Delta )^{s/2} u \big |^2 \, dx  + \frac{1}{2p +2} \int_{\R^d} |u|^{2p+2} \, dx.
\end{equation}
A natural energy space for solutions $u$ is therefore $H^s(\R^d)$, assuming that $p\in \N$ is small enough to allow for the Sobolev-imbedding $H^s(\R^d) \hookrightarrow L^{2p+2}(\R^d)$.
One should note that, due to our choice of sign in front of the nonlinearity, the energy is the sum of two non-negative terms and hence, we are in what is called the {\it defocusing regime}. 

Equation \eqref{fNLS} also enjoys the scaling invariance 
\[
u(t,x) \mapsto \frac{1}{\lambda^{s/p}}u\Big (\frac{t}{\lambda^{2s}}, \frac{x}{\lambda}\Big),
\]
which is seen to leave the homogenous Sobolev-norm $\dot {H}^\sigma(\R^d)$ invariant, provided $\sigma =\sigma_{\rm c}$, where
\begin{equation*}\label{crit1}
\sigma_{\rm c} = \Big (\frac{d}{2}- \frac{s}{p} \Big) .
\end{equation*}
One consequently refers to the Cauchy problem \eqref{fNLS} with $u_0 \in H^\sigma$ and $\sigma\ge \sigma_{\rm c}$ as $H^\sigma$ (sub-)critical, whereas the case when 
$\sigma< \sigma_{\rm c}$ is called $H^\sigma$ supercritical.
To link this notion of criticality to the energy space, one compares $\sigma_{\rm c}$ with the fractional power $s\in (0,1]$, appearing as a parameter in \eqref{fNLS}. 
Doing so implies that solutions $u(t, \cdot)\in H^s(\R^d)$ are {\it energy (sub-)critical} if $s\ge s_\ast$, where
\begin{equation}\label{crit2}
s_\ast (d,p) = \frac{d}{2} \Big (\frac{p}{p +1}\Big),
\end{equation}
while the {\it energy supercritical} regime corresponds to $s<s_\ast$. The former regime is rather well understood and the (local) 
well-posedness theory for \eqref{fNLS} in critical and sub-critical cases is, by now, broadly developed, cf. \cite{CA, HoSi, KLR}.

A systematic numerical study of equations of the form \eqref{fNLS} has first been done by the authors together with P. Markowich in \cite{KMS}. Here, we shall revisit 
these numerical simulations and expand on a particular case which has only been hinted at in \cite{KMS}. Namely, the possibility of a novel type of {\it singularity formation in the defocusing case} 
with sufficiently {\it small (fractional) dispersion}. 

This phenomenon needs to be clearly distinguished from the possibility of finite-time blow-up appearing in the case of {\it focusing} nonlinearities. 
In the focusing case, the kinetic and nonlinear potential energies in \eqref{energy} have opposite signs, which allows for the existence of solitary waves (not present in the defocusing case). 
For $s=1$, it is well known that $E(u_0)<0$ is a sufficient condition for finite-time blow-up of (strong) solutions $u\in C([0,T); H^s(\R^d))$, see, e.g., \cite{SS}. The blow-up is found to be self-similar with profiles given by the 
corresponding solitary wave solution. This basic picture has been successfully transferred to the case of fractional equations with $s\not =1$, cf. \cite{BHL} and the references given therein. 

In contrast to the focusing situation it has long been believed that defocusing equations do not suffer from any singularity formation, cf. \cite{Bo, CSS}. Indeed, it 
has recently been proved in \cite{IT1, IT2} that defocusing one-dimensional dispersive equations with {\it small} (and not necessarily localized) initial data have {\it global dispersive solutions}. 
Nevertheless, a first rigorous proof of finite-time blow-up in {\it energy supercritical defocusing NLS} ($s=1$) has been obtained by Merle et al. in their seminal paper \cite{MRRS3}.
In their work, blow-up of radial solutions $u$ is achieved by compression in the associated (classical) hydrodynamical flow of $\rho$ and $v=\nabla \varphi$, which is associated to 
$u$ via the Madelung transform $u= \sqrt{\rho} e^{i \varphi}$. By using a self-similar ansatz for both $\rho$ and $v$, and relying on analogous results for the compressible 
Euler equations \cite{MRRS1, MRRS2}, Merle et al. are able to prove the appearance of an 
oscillatory singularity with universal profile
\[
u(t,x) \Eq{|x|}0  \frac{c_1}{|x|^{{(r-1)}/{p}}} \exp\left({i \frac{c_2}{|x|^{r-2}}}\right).
\]
Here, $r=r(p,d)>2$ denotes the blow-up speed and $c_1,c_2>0$ are constants depending on $p$, $d$ and the initial data. Note that, both, the amplitude and the phase of $u$ 
become singular at the origin.

The proof of this result requires the authors of \cite{MRRS3} to work in dimensions 
$d\ge 5$ and to choose particular pairs of parameters $(d, p)$. Even though it is believed that 
the blow-up also exists for other choices of such pairs, it is unclear at this point what the range of suitable $p$ is for any given dimension $d\in \N$. Very recently, an 
extension to the case of {\it non-radial blow-up} scenarios in both $\R^d$ and $\mathbb T^d$ was provided in \cite{CLGSS} for the particular choice $(d,p)=(8,1)$. A drawback of all of these 
works is that the class of initial data for which the blow-up is proven is {\it not} explicitly specified (one only knows that there exists a finite co-dimensional manifold of smooth initial conditions). 
In particular, there is no clear threshold for $u_0$ identified beyond which blow-up will occur. This is a profound obstacle for possible numerical simulations. 

In the present paper, we shall demonstrate numerical evidence for highly oscillatory singularity formation {\it already in $d=1$}, provided the fractional (dispersive) power $s<s_\ast$ is sufficiently small. 
To this end, we shall choose several kinds of explicit initial data which are smooth, rapidly decaying and of size $\mathcal O(1)$ in ${L^\infty\cap L^2}$, cf. Section \ref{sec:crit} for details. 
From now on, we shall restrict ourselves to the situation in $d=1$ and impose
\[
u_0\in \mathcal S(\R)\subset \bigcap_{\sigma \in \R} H^\sigma(\R). 
\]
Assuming for the moment, that there exists a unique local in-time solution 
\[
\text{$u\in C([0, T); H^s(\R))$ up to some $T>0$,}
\] 
we realize that the critical Sobolev-index $s=\frac{1}{2}$ is significant in this context. 
Indeed, any solution $u(t, \cdot)\in H^s(\R)$ to \eqref{fNLS} with $s>\frac{1}{2}$, automatically satisfies $s> s_\ast$ for all $p\in \N$ and hence, we are in an energy subcritical regime. 
In addition, the conservation laws of mass \eqref{mass} and energy \eqref{energy} yield the uniform bound 
\[ 
\| u(t, \cdot) \|^2_{H^s(\R)}\le M(u_0)+ 2E(u_0), \quad \forall t\in [0,T),
\] 
which implies $T=+\infty$, i.e. global existence holds. In addition, by Sobolev imbedding, we know that if $s>\frac{1}{2}$: $H^s(\R)\hookrightarrow L^q(\R)$ for all $q\ge 2$, 
and thus all $L^q$-norms of $u(t,\cdot)$ are uniformly bounded in this case (in particular, this holds for the nonlinear potential energy).

The situation with $s<\frac{1}{2}$ is more complicated: 
in view of the Gagliardo-Nirenberg inequality, we know that for any $\theta \in (0,1)$, there exists a constant $C=C({M(u_0), E(u_0)})>0$, such that
\[
\| u \|_{L^q(\R)}\lesssim \| u\|_{L^2(\R)}^{1-\theta}  \, \| u\|^\theta_{\dot {H^s}(\R)} \lesssim C,
\]
provided $q=\frac{2}{1-2s}>2$ with $s<\frac{1}{2}$. Given a particular choice of $s<\frac{1}{2}$ this yields a uniform bound on some higher order $L^q$-norm of the solution $u(t,\cdot)$, and thus 
for the full range $q\in [2, \frac{2}{1-2s}]$. (Here, we implicitly assume that the nonlinear power $p\in \mathbb N$ is 
such that $q=2p + 2$ is within this admissible range, in order to guarantee that the nonlinear potential energy of $u$ is well-defined.) 
Clearly, this higher order $L^q$-bound deteriorates as $s\to 0$ and eventually only yields the obvious bound on the $L^2$-norm (or mass) of $u$. 
For small enough $s<\frac{1}{2}$ this consequently leaves open the possibility 
that the solution might develop a singularity in some higher $L^q$-norm where $2<q\le \infty$. 

Let us note that $s<\frac{1}{2}$ does not necessarily correspond to 
an energy supercritical regime. Indeed, the smallest energy critical index in $d=1$ is found to be $s_\ast = \frac{1}{4}$, obtained in the case of a cubic nonlinearity (where $p=1$). More generally, 
we have the following range of critical indices
\[
s_\ast =\frac{p}{2p+2} \in \Big [\frac{1}{4}, \frac{1}{2}\Big) \quad \text{in $d=1$.} 
\]

\medskip

This paper is now organized as follows: in Section \ref{sec:num} we introduce a semiclassical rescaling for fNLS and 
review the numerical methods deployed in our simulations. In Section \ref{sec:linear} we discuss the linear fractional Schrödinger equation and 
some of its basic dynamical properties. Numerical simulations 
for various initial data in the energy (sub-)critical regime are presented in Section \ref{sec:crit}. Finally, in Section \ref{sec:supcrit} we shall provide numerical 
evidence for a singularity formation in various energy supercritical cases. In the appendix, we collect some remarks on the comparison of our numerical scheme with other numerical methods.


\section{Semiclassical rescaling and numerical methods}\label{sec:num}

\subsection{Semiclassical rescaling} 

Recall that the blow-up phenomenon in \cite{MRRS3} is given by a highly oscillatory singularity obtained from the underlying classical hydrodynamical flow. 
Motivated by these ideas, we shall consider \eqref{fNLS} in a {\it semi-classical regime} which imprints high-frequency oscillations onto the solution. 

To this end, let $\eps \sim \hbar \ll 1$ be a small (dimensionless) {\it semi-classical parameter}, and consider slowly varying initial data of the form $$u_0(x) = \upsilon(\eps x),$$ 
where $\upsilon \in \mathcal S(\R)$ is some given $\eps$-independent profile of size $\mathcal O(1)$. 
As $\eps \to 0$ the initial data approach the constant value $\upsilon(0)$. Hence, in order to see nontrivial effects within the solution, we will need to consider macroscopic length- and time-scales 
of order $\mathcal O(1/\eps)$. As in \cite{KMS} we rescale $x\mapsto  \tilde x= \eps x$, $t\mapsto  \tilde t = \eps t$, and denote the new unknown by
\[
\tilde u(\tilde t,\tilde x) = u \Big(\frac{\tilde t}{\eps}, \frac{\tilde x}{ \eps}\Big).
\] 
The latter solves a semi-classically scaled equation, where we discard the $\tilde{ }$ again, for simplicity:
\begin{equation}\label{resfNLS}
i \eps \partial_t u  = (- \eps^2 \Delta )^s u + |u|^{2p} u, \quad u_{\vert{t=0}} = \upsilon(x),
\end{equation}
At least for short times $0<|t|<t_0$, we expect the solution to exhibit a highly oscillatory behavior of the form
\begin{equation}\label{wkb}
u(t,x) \Eq\eps0 a(t,x)e^{i \varphi(t,x)/\eps},
\end{equation}
where $\varphi$ is some real-valued phase and $a$ denotes a slowly varying amplitude. 

\begin{remark}
For $s=1$ it is well known that an ansatz of the form \eqref{wkb} allows one to study the semi-classical asymptotics of $u$ as $\eps \to 0$ via an 
Euler-type system obtained for the density $\rho = |a|^2$ and the velocity $v= \nabla \varphi$. The latter is known 
to be equivalent to the original Schr\"odinger equation, cf. \cite{Car}. In the case $s\not =1$ no such equivalent formulation (or Madelung transform) is known at this point. 
\end{remark}

A particular situation appears in a regime with very small $s\ll 1$, which can be considered as a {\it small dispersion limit} for \eqref{resfNLS}. Indeed, we see that (at least formally)
\[
(-\eps^2 \Delta)^{s} u \to {\bf 1}u , \quad \text{as $s\to 0$,}
\]
resulting in the following limiting equation
\begin{equation}\label{limeq}
i\eps \partial_t u_{\rm lim} = (1 + |u_{\rm lim}|^{2p}) u_{\rm lim}, \quad {u_{\rm lim}}_{\vert{t=0}} = \upsilon(x).
\end{equation}
The latter is {\it non-dispersive} and its solution is explicitly given by
\[
u_{\rm lim}(t,x) = \upsilon(x) e^{-it (1+|\upsilon(x)|^{2p})/\eps}.
\]
This implies that the particle (or energy) density $|u_{\rm 
lim}(t,x)|^2 = |\upsilon(x)|^2$ for all $t\in \R$ and is hence, $\eps$-independent, while the solution itself is seen to exhibit
(nonlinear) oscillation with frequency $\nu = \mathcal O(\tfrac{1}{\eps})$.

\begin{remark} In our numerical simulations below, we shall typically choose $\eps = 0.1$. This allows us to accurately compute up to times of order $T\simeq 20$. 
\end{remark}


\subsection{Numerical method for computing the time-evolution}\label{sec:numtime}

To solve equation (\ref{resfNLS}) numerically, we shall use a 
Fourier-spectral method similar to the one in \cite{KMS}. This choice 
is motivated by the highly oscillatory nature of our problem with 
steep gradients and thus the need for very high numerical resolution. We refer to \cite{KSb}
for other possible approaches to numerically simulate fractional 
derivatives in less oscillatory situations.

To be more precise, we shall construct the numerical solution to \eqref{resfNLS} as $$u(t, x) = (\mathcal F^{-1}\widehat u)(t, x),$$ where $\widehat u$ solves
\begin{equation}
	i \eps \partial_{t}\widehat{u}=\eps^{2s}|\xi |^{2s}\widehat{u}+\mathcal F({|u|^{2}u}),\quad \widehat u_{\vert{t=0}} = \widehat \upsilon(\xi)
	\label{fNLSfourier}.
\end{equation}
We thereby consider $x\in  [-L\pi,L\pi]$ with $L\gg 1$,
and impose {\it periodic boundary conditions} to approximate $\widehat{u}(\cdot, \xi)\approx \widehat u(\cdot, k)$, i.e. the discrete Fourier-transform of $u(\cdot, x)$ 
with frequency-variable $k \in \frac{1}{2L}\Z$ (by a slight abuse of notation, we shall denote both the Fourier transform and its discrete version by
the same symbol). 

The function $\widehat u(\cdot ,k)$ can be conveniently computed using the Fast Fourier transform (FFT). We thereby apply 
the standard discretization for the FFT, using
\[
x_{n}=\big(-\pi+\tfrac{2n\pi}{N}\big)L, \quad  n=1,\dots,N\in \N,
\] 
with the corresponding discrete Fourier variable
\[
k=\tfrac{1}{L}\big(-\tfrac{N}{2} +1,\dots,\tfrac{N}{2}\big)\in \tfrac{1}{2L}\Z.
\]
For functions $ u(\cdot, x)$ that are analytic in $x\in  [-L\pi,L\pi]$ within 
numerical precision (here, roughly $10^{-16}$), it is well known that the numerical error 
decreases exponentially with $N$, see \cite{trefethen}. However, it was shown in \cite{KSb, KS} that this exponential decrease in $N$ cannot be expected to hold for fractional powers $s\notin \mathbb{N}$, since the Fourier-symbol 
$|k|^{2s}$ is no longer regular. Indeed, a Paley-Wiener argument implies that the inverse Fourier transform of 
$|k|^{2s}\widehat{u}(\cdot, k)$ is not rapidly decreasing, even if $\widehat{u}(\cdot, k)$ is in 
the class of Schwartz functions w.r.t. $k$. In turn this  
means that machine precision {\it cannot} be reached for non-integer $s\in (0,1]$, see \cite{KS}. 
Indeed, in all of our simulations the numerical error remains of 
the order of $10^{-5}$, even when using a very high spatial resolution. This (relatively large) numerical error is nevertheless acceptable for our purposes. In addition, the use of Fourier-methods 
allows us to implement
a numerical algorithm for singularity tracing in the frequency domain, see Section \ref{sec:singtr} below.

The FFT discretization of \eqref{fNLSfourier} leads to an $N$-dimensional system of ordinary 
differential equations which for $s\simeq 1$ are \emph{stiff} (meaning that explicit methods are usually inefficient in view of 
the required stability conditions).  
However, since the stiff terms are in the diagonal part of the evolutionary $N\times N$-matrix, there are many efficient 
integrators for such systems, cf. \cite{etna}. In the present work, 
we are mainly interested in values of $s\le \frac{1}{2}$, and thus {\it small to moderate stiffness}, which 
makes it possible to employ the classical {\it explicit fourth order 
Runge-Kutta method} (RK4). The latter is three to four times 
faster than the implicit code (IRK4) used in our previous work \cite{KMS} and yet, all results 
obtained in the present paper agree with the ones in \cite{KMS} up to errors of order $10^{-3}$, or better.


The accuracy of the time-integration is 
controlled via the conservation of the numerically computed energy 
\eqref{energy} and mass \eqref{mass}. Due to unavoidable numerical errors, the numerical energy $E_{\rm num}(t)\equiv  E((\mathcal F^{-1}\widehat u)(t, \cdot))$ will depend on time. 
In our simulations, the computed relative energy of the solution
\begin{equation}
    \Delta_E=\Big| \frac{E_{\rm num}(t)}{E_{\rm num}(0)}-1\Big|
    \label{Delta}
\end{equation}
remains of the order $10^{-7}$ at worst, while the mass $M_{\rm 
num}(t)$ is conserved up to errors of order $10^{-10}$. In addition, we redo all of our critical numerical experiments using 
double the numerical resolution in the second run and confirming that 
the results agree within the desired accuracy of order $10^{-3}$. Comparison of our RK4 scheme with other fourth-order methods reaffirm this accuracy (see the discussion in the appendix). 

\begin{remark} In \cite{etna} it was argued that the computed energy 
conservation overestimates 
the accuracy of the numerical solution by one to two orders of magnitude. 
\end{remark}


\subsection{Singularity tracing in the complex plane}\label{sec:singtr}

An advantage of Fourier methods is that they allow us to trace 
singularities in the complex plane via the asymptotic behavior of the 
Fourier transform. This was first used in a numerical context by 
Sulem et al. in \cite{SSF} and has later been successfully applied in \cite{KR2013a, KR2013b}. 

Assume that $u(\cdot, x)$ can be analytically extended to a function $u(\cdot, z)$ within a strip around the real axis inside the complex plane. 
Suppose further that at a certain time $t\in \R$ this extended 
function $u$ has a single essential singularity at $z_0(t)\in \C$, such that 
\[
\text{$u(t, z) \sim (z-z_{0}(t))^{\mu(t)} \ $ for $z\sim z_{0}(t)$, with $\mu(t) \not \in \mathbb N$.} 
\]
Classical results in asymptotic analysis (see, e.g, \cite{asymbook}) then imply that 
the Fourier transform has the following asymptotic behavior, as $|\xi |\to\infty$:
\begin{equation}
    |\widehat{u}(t, \xi)|\sim 
    \frac{1}{|\xi |^{1+\mu(t)}} e^{-\delta(t) |\xi |},\quad \text{where $\delta(t)=\text{Im}\, z_{0}(t)$.}
    \label{fourasymp}
\end{equation}
In \cite{KR2013a, KR2013b} the behavior of the FFT 
coefficients $\xi = k \in \frac{1}{2L}\Z$ within the numerical solution $u (t, \cdot)$ (which asymptotically should be of the form \eqref{fourasymp})
was used to quantitatively identify the time $t=t_{\rm f}$ where the 
singularity $z_0(t)=x_0(t)+i \delta(t)$ hits the real axis, i.e., where $\delta(t_{\rm f}) = 0$ and $u(t_{\rm f}, x_0)$ becomes singular. The quantity $\delta\equiv \delta(t_{\rm f})$ is thereby identified 
by fitting the logarithm of the modulus of the FFT coefficients  via linear regression. 
Since $|\xi|^{\mu+1}$ is only a logarithmic correction to the exponential behavior in \eqref{fourasymp}, the quantity $\mu\equiv \mu (t_{\rm f})\in \R$ which characterizes the type of the singularity, is 
unfortunately less reliably determined than $\delta$. However, in the 
case of the defocusing NLS, it was successfully used to detect the well-known cusp singularity $\mu = \frac{1}{3}$ 
in semiclassical limits, cf. \cite{KR2013b}.

In practice, the minimal resolved 
distance via Fourier methods is given by
\begin{equation*}
   m:=2\pi \frac{L}{N}\ll 1,
\end{equation*}
with $N\in \N$ being the number of 
FFT coefficients. This means that once $\delta (t) \leq  m$, it is compatible 
to zero with numerical accuracy $\delta \simeq 0$. It was shown in \cite{KR2013b} that 
the best results are obtained when the numerical code is stopped at a 
final time $t=t_{\rm f}$ at which $\delta (t)$ actually vanishes, 
and this is the criterion to be applied also in the present paper. In 
our numerical simulations, we always find $0< t_{\rm f} \le 20$. 

\begin{remark} Note that the expression (\ref{fourasymp}) only holds 
in the presence of a single such singularity. If there are several 
singularities $z_{j}(t)\in \mathbb{C}$, $j=1,\ldots,J$, 
each of them contributes asymptotically, resulting in
\[
\widehat u (t,z) \sim \sum_{j=1}^{J} \frac{e^{z_j (t) |\xi|}}{|\xi |^{1+\mu_j(t)}} , \quad \text{as $|\xi |\to\infty$.}
\]
This leads to interferences due to the presence of oscillations whenever $\text{Re}\, z_j\not = 0$. In such 
cases the vanishing of the first parameter $\delta_{j}$ will still 
indicate that a singularity has hit the real axis, but it will be 
difficult to identify the corresponding $\mu_{j}$. 
\end{remark}


\section{Linear fractional Schr\"odinger equation}\label{sec:linear}
In this section we shall study the linear case, i.e.,
\begin{equation}\label{fschr}
i \eps \partial_t u = (- \eps^{2} \Delta )^s u, \quad u_{\vert{t=0}} = \upsilon,
\end{equation}
in order to gain some basic insight on the influence of the fractional dispersion within the time evolution. The solution to (\ref{fschr}) is given by
\begin{equation}
	u (t, x)= \mathcal F^{-1} \left( \widehat{\upsilon}(\xi)e^{-i\eps^{2s-1}|\xi|^{2s}t}\right)
	\label{fschrfourier}.
\end{equation}
No time integration scheme is needed in this case; the inverse Fourier transform of 
\eqref{fschrfourier} is just approximated for any given time $t$ by an 
inverse FFT. For $\eps =0.1$, we will do so, by using $N=2^{14}$ FFT 
modes for $x\in [-10 \pi,10 \pi]$. As initial data we consider here
\begin{equation}\label{inilin}
\upsilon(x)={\mbox{sech}(x)},
\end{equation} 
but the pictures remain qualitatively similar for other choices of $\upsilon$ (such as those given in Section \ref{sec:crit} below). 

\begin{figure}[htb!]
 \includegraphics[width=0.49\textwidth]{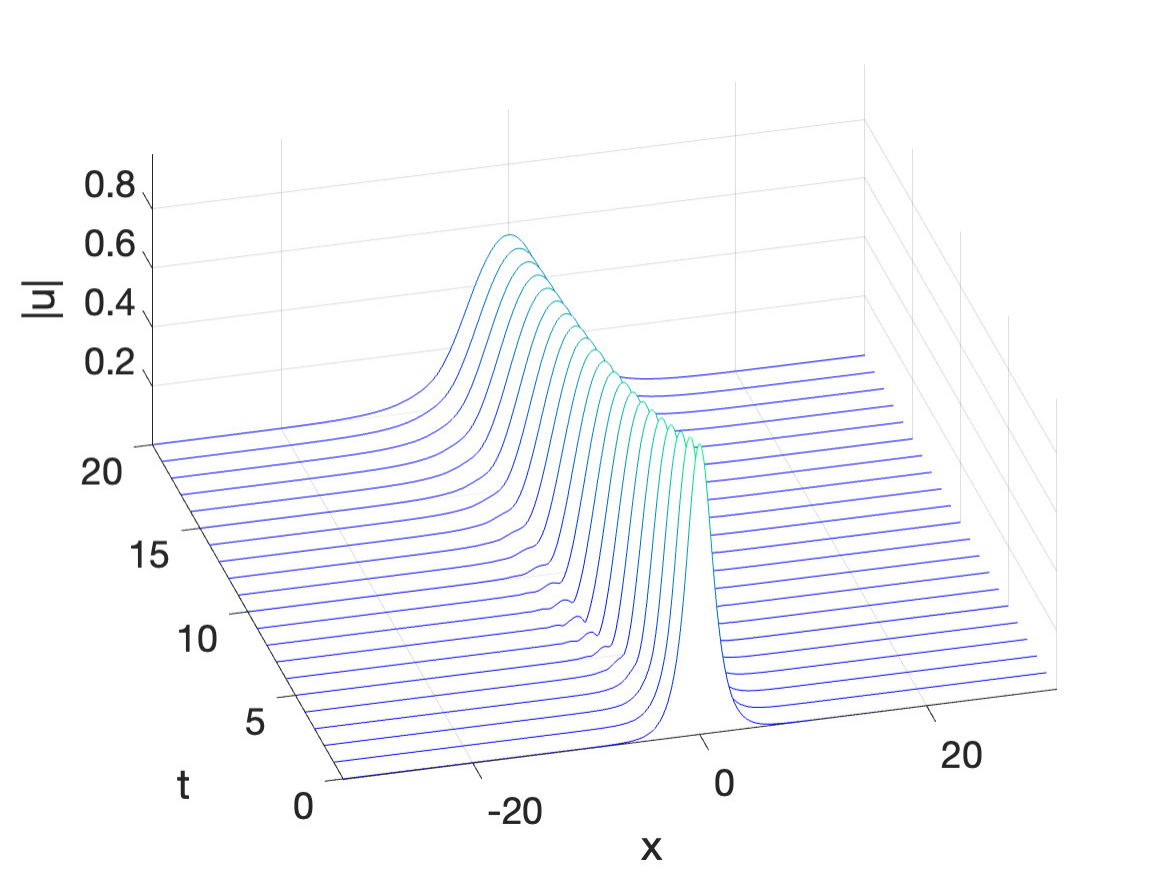}
 \includegraphics[width=0.49\textwidth]{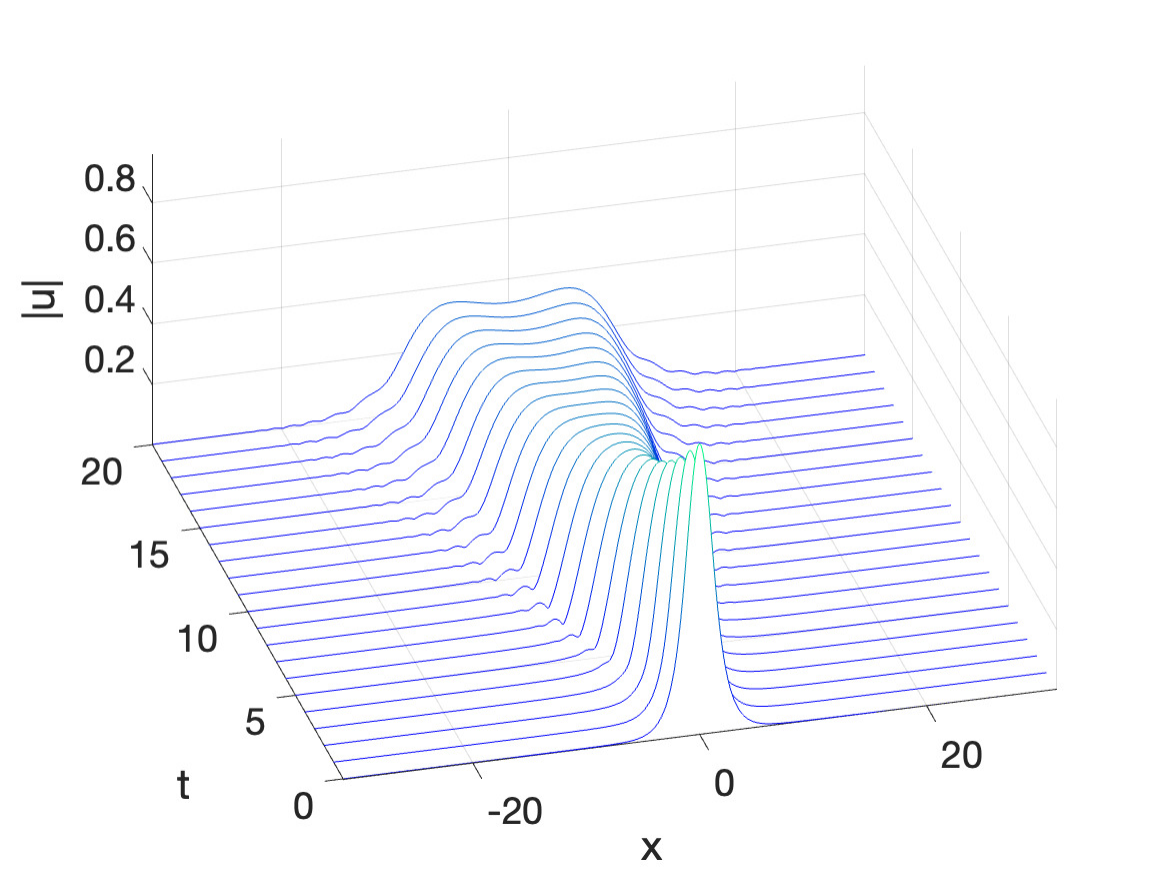}
 \caption{Solution to the linear equation 
 \eqref{fschr} with $\eps=0.1$ and initial data \eqref{inilin}: on the left for $s=1$, on the right for $s=0.8$.  }
 \label{fschr1}
\end{figure}

First, we compare the solution to the classical Schr\"odinger 
equation with $s=1$ to the case with a slightly smaller fractional power 
$s=0.8$. In Fig.~\ref{fschr1}, the modulus $|u|$ of both solutions 
is shown as a function of $t$ and $x$. 
On the left, where $s=1$, the 
initial hump is seen to be slowly dispersed as $t$ becomes large. On the right, it can 
be seen that this effect is enhanced for $s=0.8$. 

\begin{figure}[htb!]
 \includegraphics[width=0.49\textwidth]{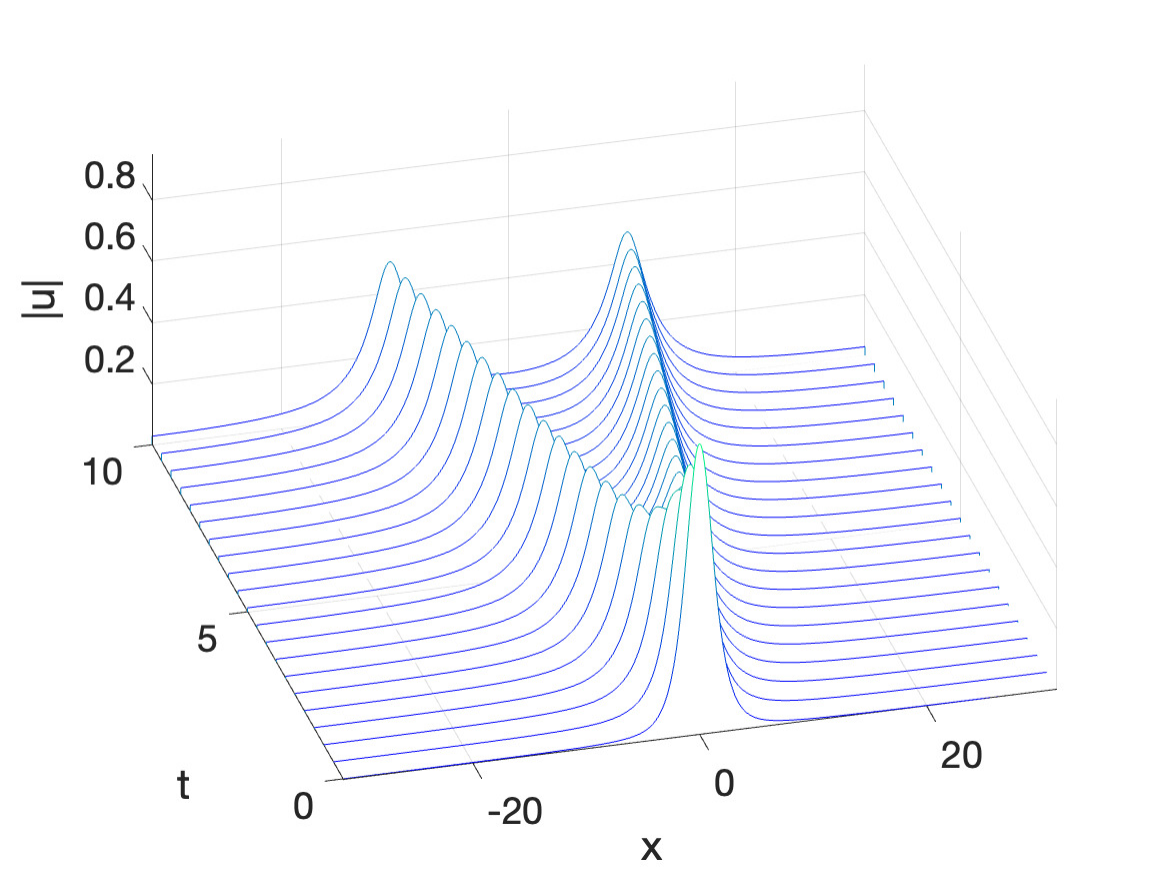}
 \includegraphics[width=0.49\textwidth]{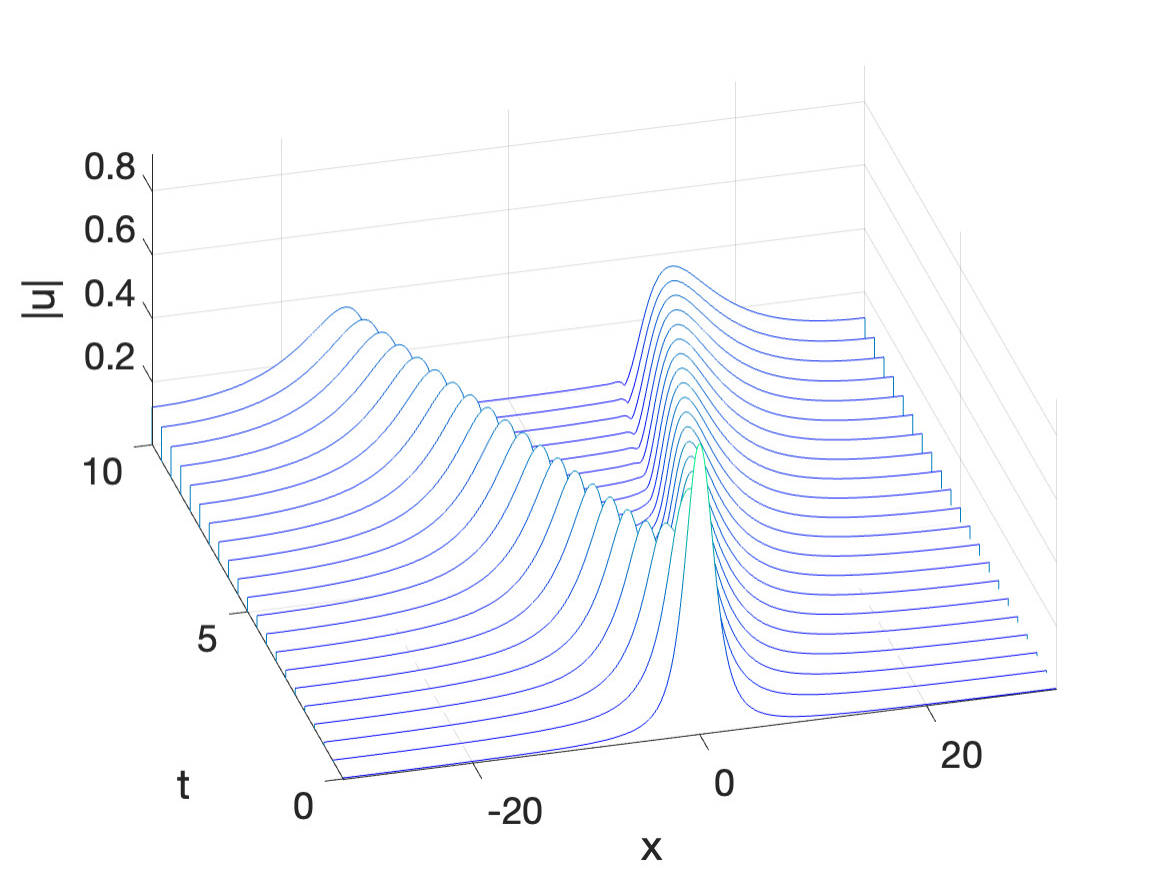}
 \caption{Solution to the linear equation 
 \eqref{fschr} with $\eps=0.1$ and initial data \eqref{inilin}: on the left for $s=\frac12$, on the right for $s=0.4$.  }
 \label{fschr05}
\end{figure}

For smaller values of $s\in (0,1]$ it can be observed that the initial hump splits 
into two humps traveling towards $x\to \pm\infty$. In Fig.~\ref{fschr05}, we show the case $s=\frac12$ and $s= 0.4$, respectively. 
The former corresponds to a {\it hyperbolic} Schr\"odinger equation with dispersion relation $\lambda(\xi)=|\xi|$, i.e.
\begin{equation}\label{hyper}
i \partial_t u = \sqrt{-\Delta} u , \quad u_{\vert{t=0}} = \upsilon.
\end{equation}
Note that in this case, the small semiclassical parameter $\eps$ cancels out of the equation, while for $s < \frac{1}{2}$ the power of $\eps$ appearing 
in formula \eqref{fschrfourier} becomes negative. For $s=0.4$ the solution for $t>1$ appears to be almost separated into two humps with 
rather steep interior fronts (cf. the right picture in Fig.~\ref{fschr05}). 
Since the equation is linear, the solution nevertheless stays smooth at all times. 

\begin{figure}[htb!]
 \includegraphics[width=0.49\textwidth]{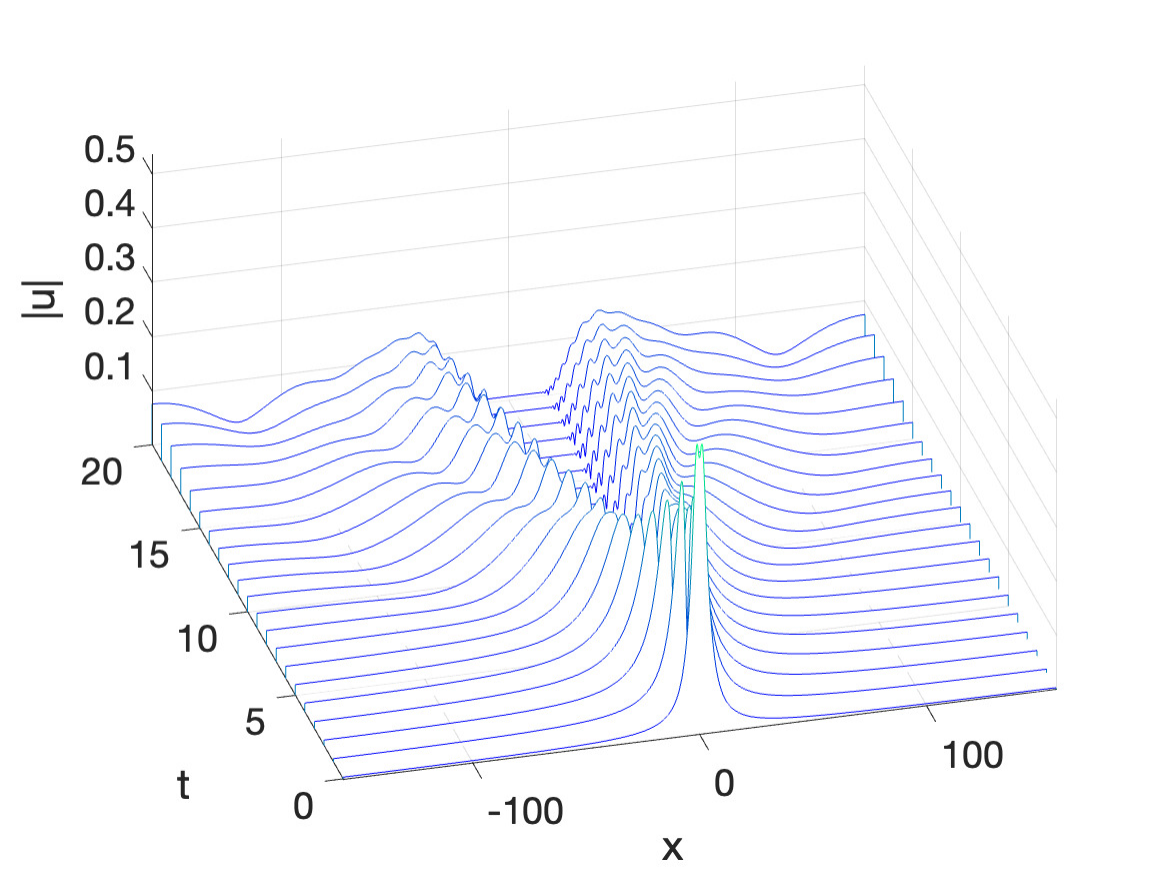}
 \includegraphics[width=0.49\textwidth]{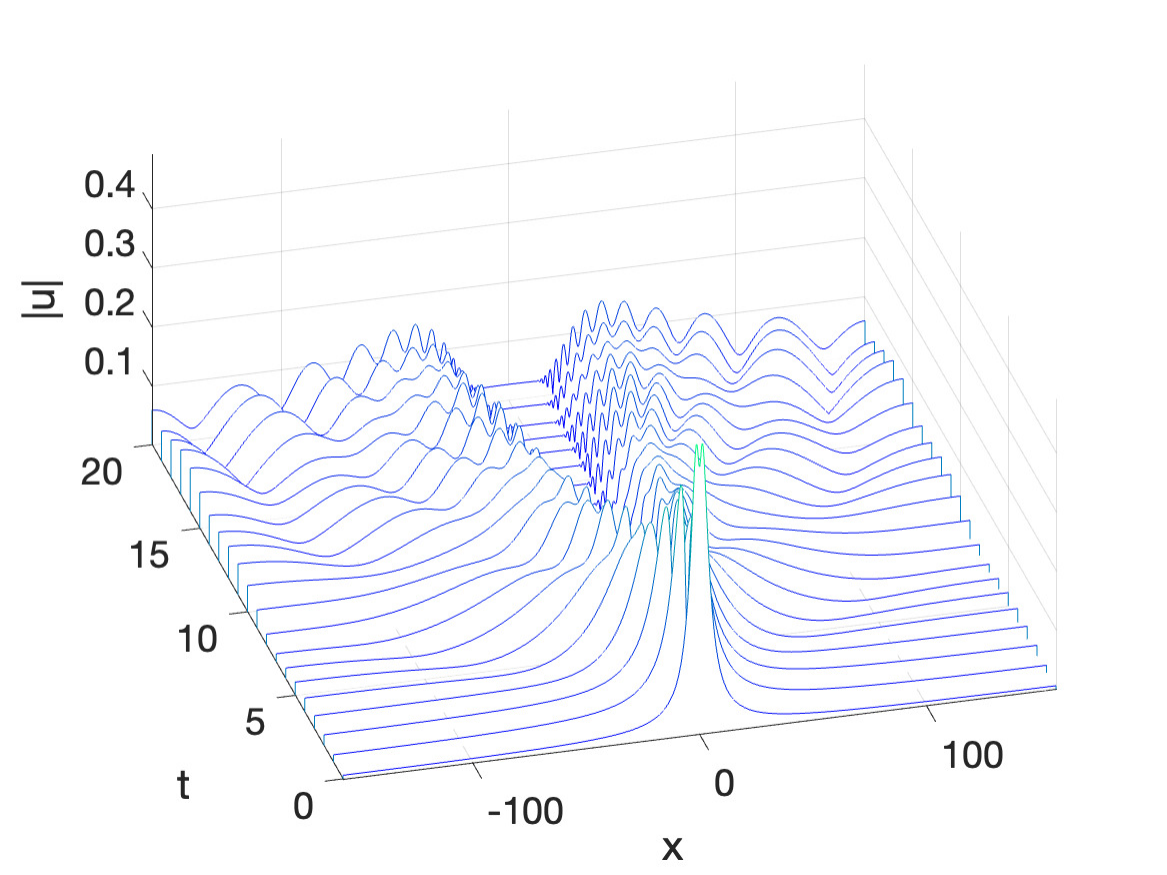}
 \caption{Solution to the linear equation 
 \eqref{fschr} with $\eps=0.1$ and initial data \eqref{inilin}: on the left for $s=\frac14$, on the right for $s=\frac15$.  }
 \label{fschr025}
\end{figure}

For even smaller values of the dispersion, the solution still 
develops two main humps, but the latter now decompose into several 
smaller crests. 
In Fig.~\ref{fschr025}, we show on the left the case with $s=\frac{1}{4}$, 
which corresponds to the energy critical regime for the corresponding 
nonlinear model \eqref{resfNLS} with $p=1$. On the right of the same figure we show the case with 
$s=\frac15$. In this particular situation, the dispersion appears to lead to more than 
the two humps since there is a more agitated wave pattern to be 
observed, im particular at the interior boundary of the humps. 

These pictures can be partly explained by the fact that we expect the bulk of the particle density $\rho(t,x)=|u(t,x)|^2\sim |a(t,x)|^2$ 
to approximately travel along the rays obtained via geometric optics, i.e.
\[
x\mapsto X(t; x) = x +t \nabla \lambda(\xi)
\]
The slope of these rays is given by the group velocity $v(\xi) = \nabla \lambda(\xi)$, which in $d=1$ is equal to
\[
v(\xi) =  2s\,  {\rm sgn}(\xi) |\xi|^{2s-1 }.
\]
In particular, in the case of the hyperbolic equation \eqref{hyper} with $s=\frac{1}{2}$, the group velocity is $v(\xi) = {\rm sgn}(\xi) = \pm 1$, which explains the 
splitting of the initial hump in Fig. \ref{fschr05}. For $s>\frac{1}{2}$ the solution is seen to interpolate between these two main directions of propagation in a smooth way. For $s<\frac{1}{2}$, however, 
the group velocity becomes singular at the origin. In this regime the effect of dispersion is ``inverted" and wave packets with frequencies $|\xi|<1$ will 
move faster than those with frequencies $|\xi|>1$. We therefore expect a type of compression effect within the solution, as the slow frequencies within $u$ try to ``overtake" the high frequencies. 


\section{Solution of the fractional NLS: the (sub-)critical regime}\label{sec:crit}

In this section we shall expand on some of the numerical experiments from \cite{KMS} and consider specifically the fNLS with cubic nonlinearity ($p=1$), i.e.
\begin{equation}\label{critfNLS}
i \eps \partial_t u  = (- \eps^2 \Delta )^s u + |u|^{2} u, \quad u_{\vert{t=0}} = \upsilon,
\end{equation}
for which we shall study the behavior of its solution in energy subcritical and critical cases. 
For $s=1$ and $\eps \ll 1$ this corresponds to the well known semiclassical limiting regime for 
cubic nonlinear Schr\"odinger equations, a particular case of the small dispersion limits studied in \cite{DGK2}. 


\subsection{The subcritical regime} We first consider the subcritical regime with $s>\frac{1}{2}$, for which we use $N=2^{14}$ FFT 
modes to resolve $x\in [-10 \pi,10 \pi]$, and integrate equation \eqref{resfNLS} using $N_{t}=10^{3}$ time steps. For the initial data $\upsilon$ 
we shall consider three particular cases of smooth and rapidly decreasing functions. Namely,
\begin{equation}\label{ininonlin}
{\rm (a)} \ \upsilon(x) = {\mbox{sech}(x)}\qquad {\rm (b)}\  \upsilon(x) = e^{-x^2}\qquad {\rm (c)}\  \upsilon(x) = e^{-x^4}.
\end{equation}

\begin{figure}[htb!]
 \includegraphics[width=0.49\textwidth]{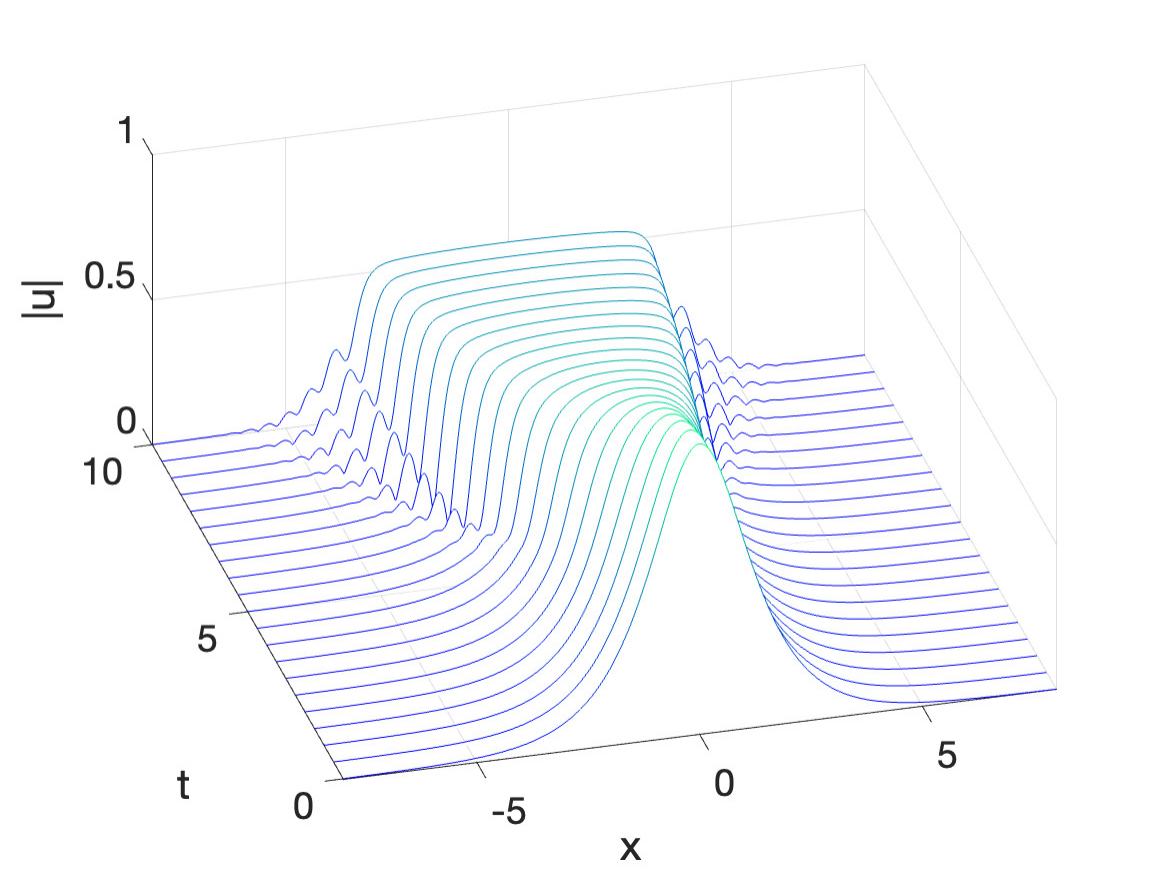}
 \includegraphics[width=0.49\textwidth]{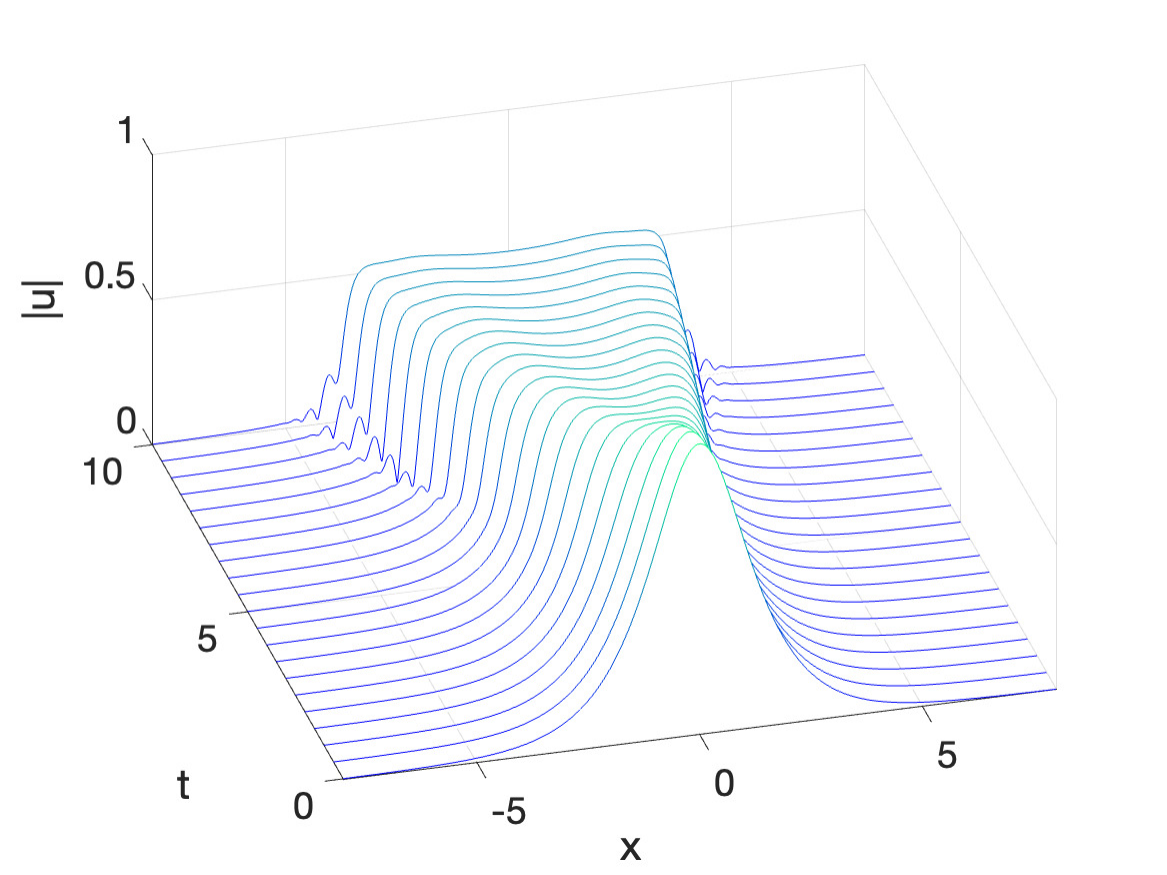}
 \caption{Solution to the cubic fNLS \eqref{critfNLS} with $\eps=0.1$ and initial data \eqref{ininonlin}(a): on the left for $s=1$, on the right for $s=0.8$.  }
 \label{fNLSs1}
\end{figure}

In the left picture of Fig.~\ref{fNLSs1}, we show the modulus of the solution $u$ to the semiclassically scaled NLS ($s=1$) with $\eps=0.1$ and initial data \eqref{ininonlin}(a) as a function of $t$ and $x$. 

As in the linear case, the dispersion leads 
to a broadening of the initial profile, but the nonlinearity implies a 
steepening of the gradients at the outer fronts. The small dispersion coefficient $\eps =0.1$ 
leads to oscillations near these fronts. 
As can be seen on the right of the 
same figure, the situation for $s=0.8$ is qualitatively similar. In both of these cases, the general 
behavior can be interpreted as a combination of the linear dynamics (cf. Fig.~\ref{fschr1}) with additional small 
oscillations in the region of strong gradients.   


\subsection{The energy critical regime: emergence of turbulence}

For computing the solution in the energy critical case $s=\frac14$, we use $N=2^{19}$ FFT modes for $x\in 
[-20\pi,20 \pi]$ and $N_{t}=10^{5}$ time steps. 

\begin{figure}[htb!]
 \includegraphics[width=0.8\textwidth]{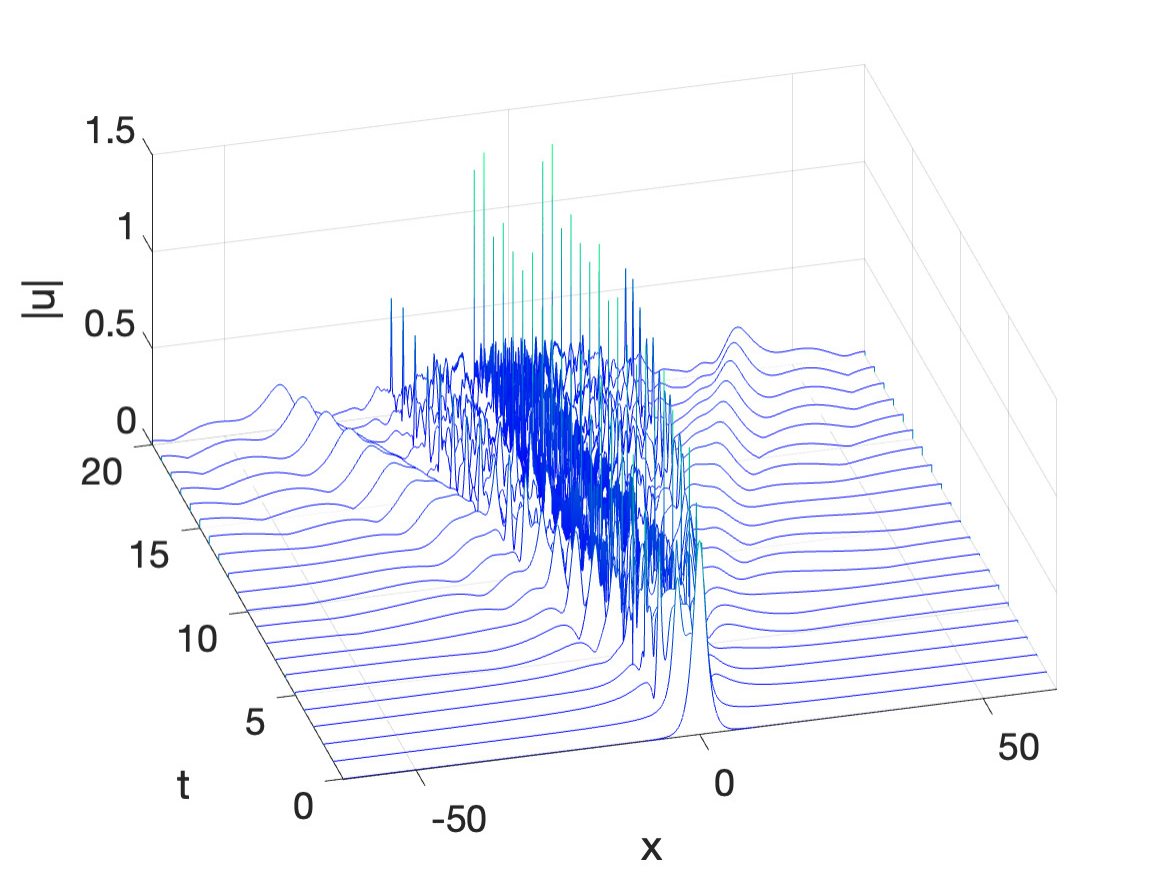}
 \caption{Solution to the energy-critical cubic fNLS \eqref{critfNLS} with $s=\frac{1}{4}$, $\eps=0.1$, and initial data \eqref{ininonlin}(a).}
 \label{fNLSs025water}
\end{figure}

In Fig.~\ref{fNLSs025water} the modulus of the solution with initial data \eqref{ininonlin}(a) is shown as a function of $t$ and $x$. It can be seen that $|u|$ exhibits rapid 
oscillations in between some very thin peaks. 
This behavior is consistent with the onset of a phenomenon called {\it weak turbulence}, i.e. a low-to-high frequency cascade within $\widehat u(t, \xi)$, in which the mass $M(u)$ and energy $E(u)$  
remain conserved, but $\widehat u(t, \xi)$ shifts to increasingly higher values in $|\xi|$, or equivalently, to finer spatial scales in $x$, as $t\to +\infty$. 
\begin{remark} 
For cubically nonlinear Schr\"odinger equations ($s=1$), the regime of weak turbulence has been rigorously studied in recent years, cf. \cite{BGHS, CSTT, Ha, Ku} 
and the references therein. 
In accordance with these studies, turbulent solutions $u$ are to exhibit supercritical $H^s$-norms growth as $t\to \infty$, cf. Fig 8 and 9 below. 
\end{remark}

\begin{figure}[htb!]
 \includegraphics[width=0.49\textwidth]{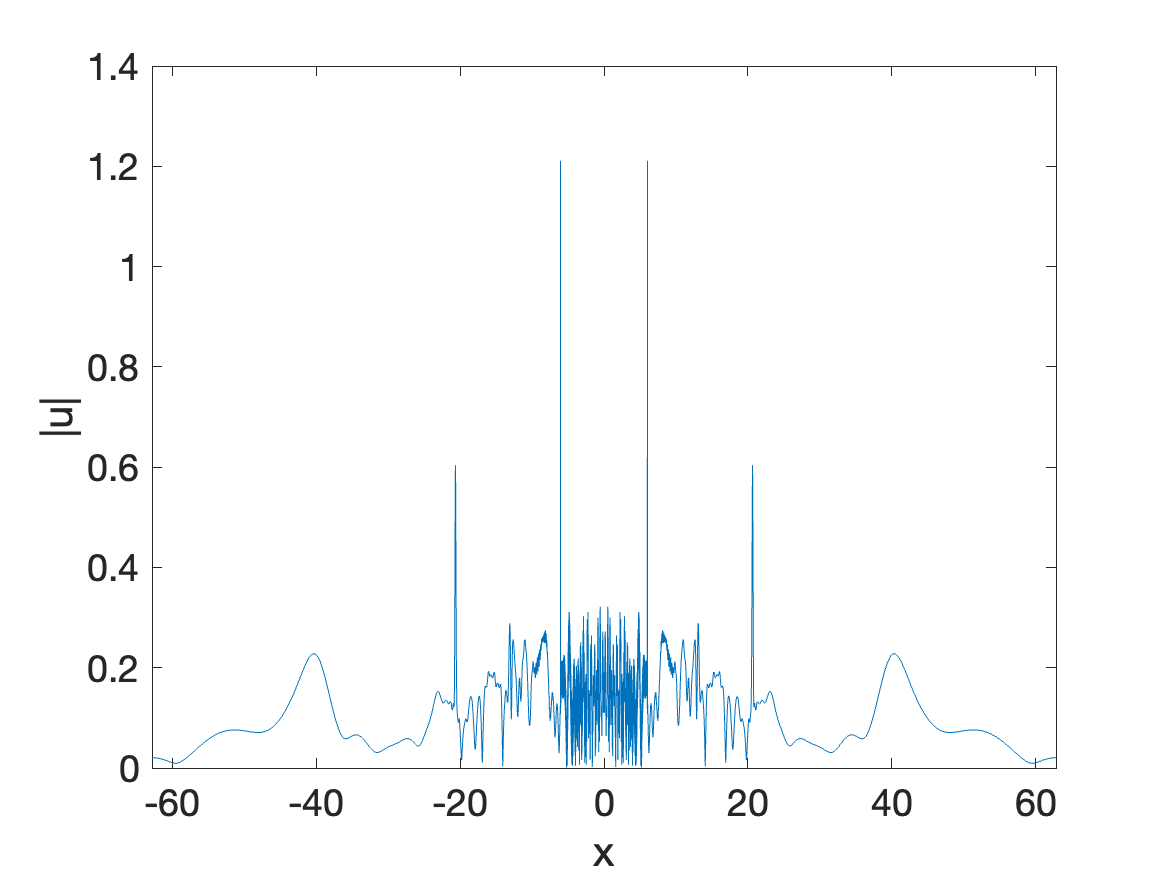}
 \includegraphics[width=0.49\textwidth]{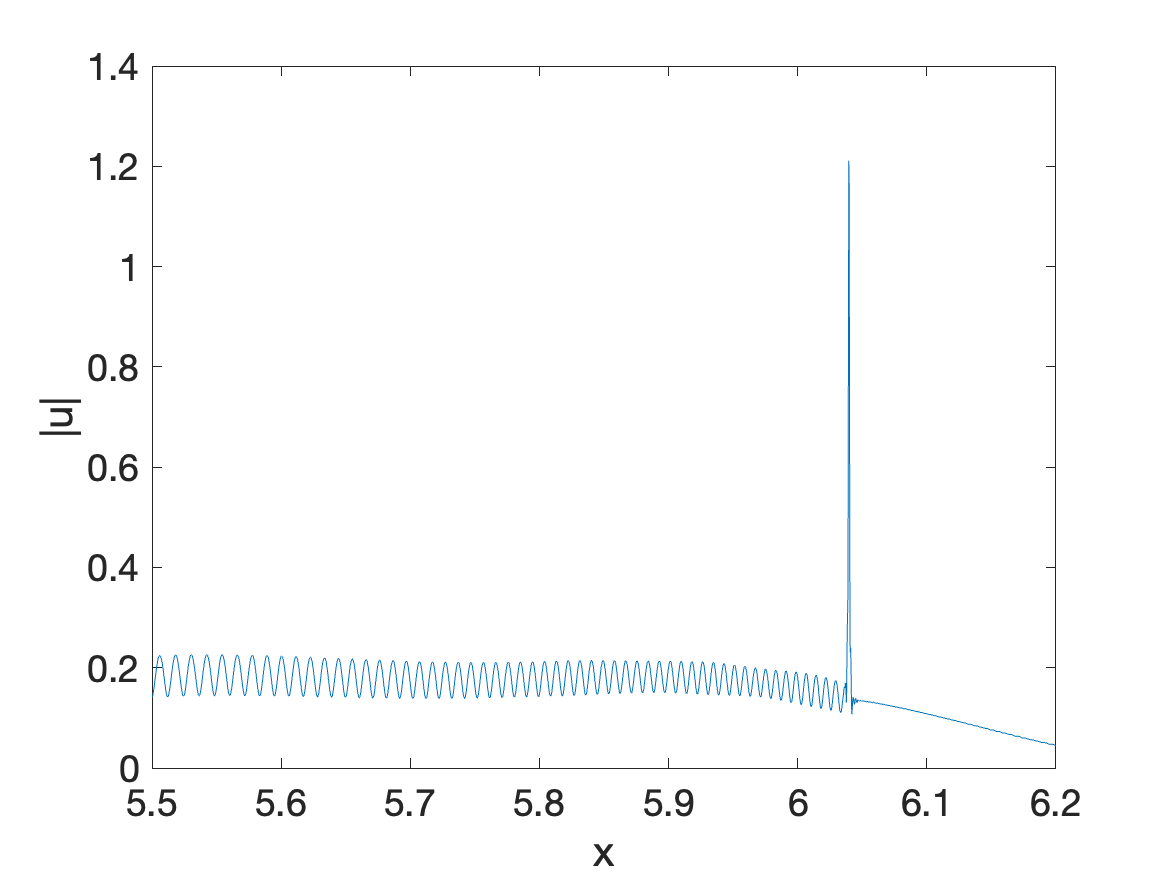}
 \caption{Left: Solution to the cubic fNLS \eqref{critfNLS} with $s=\frac{1}{4}$, $\eps=0.1$, and initial data \eqref{ininonlin}(a) at $t=20$. Right: a close-up 
 near one of the peaks.}
 \label{fNLSs025t20}
\end{figure}

A plot of the (modulus of the) solution of Fig.~\ref{fNLSs025water} at time $t=20$ is 
given in Fig.~\ref{fNLSs025t20}. It shows two large narrow peaks with finite amplitude and a zone of rapid 
oscillations in between. On the right of the same figure, a close-up 
of one of the peaks is shown, confirming that the solution is still well resolved numerically. We can see that the solution is still 
regular at this time. In particular, the oscillations do not 
increase their frequency near the peak. 

\begin{figure}[htb!]
 \includegraphics[width=0.49\textwidth]{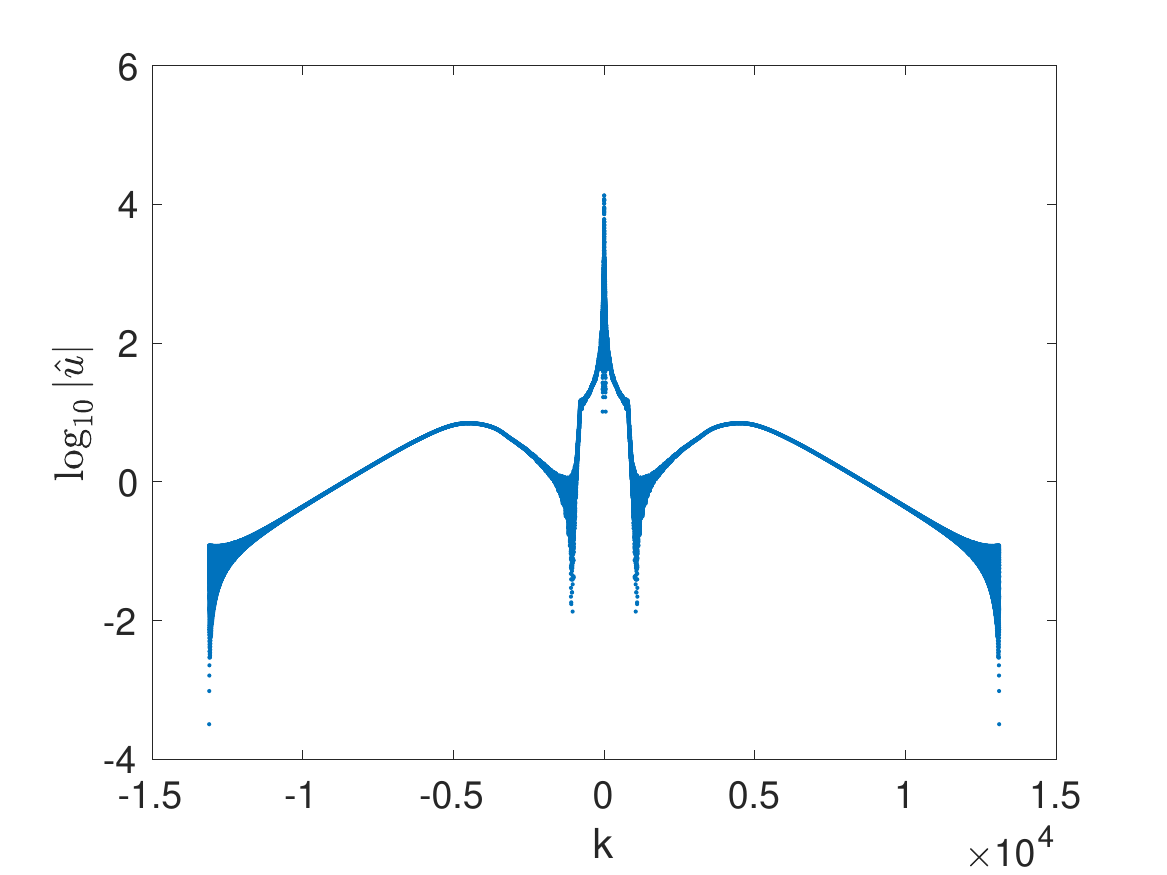}
 \includegraphics[width=0.49\textwidth]{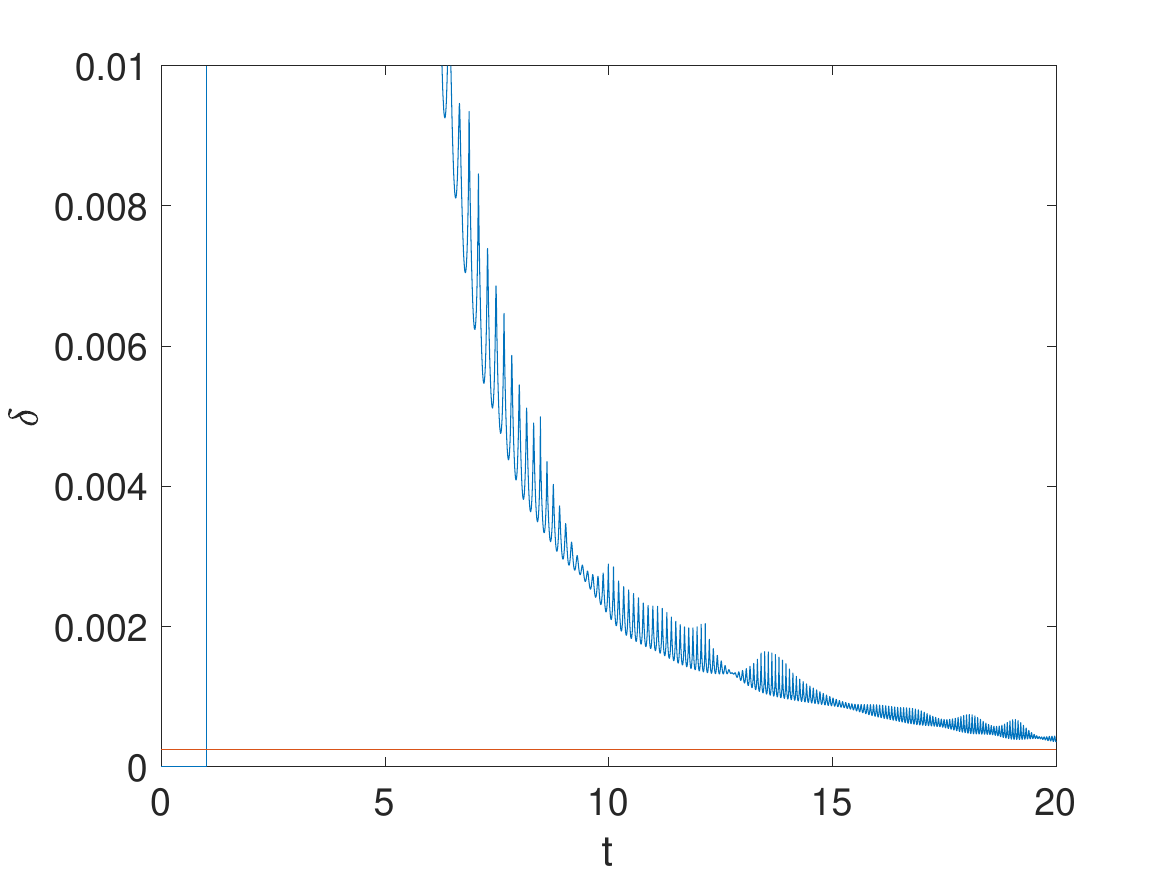}
 \caption{Left: FFT coefficients of the solution from Figure \ref{fNLSs025t20}. 
 Right: the fitting parameter $\delta$ as a function of time.}
 \label{fNLSs025t20fourier}
\end{figure}

The regularity of $u$ is confirmed by the fact that the 
FFT coefficients of $\widehat u(t,k)$ at $t=20$ span about five orders of magnitude, cf. the 
left of Fig.~\ref{fNLSs025t20fourier}. The fitting parameter $\delta(t)$, appearing in the asymptotic formula  
\eqref{fourasymp}, is shown as a function of time on the right of the same figure. 
We see that $\delta(t)$ is clearly decreasing, indicating that an essential singularity within $u(t,z)$ approaches the real axis. 
However, at $t=20$ the parameter $\delta > m$ and hence, still above the minimal 
distance $m$ resolved by the FFT discretization. 

\begin{remark} Numerically it cannot be excluded that $\delta(t) $ will eventually drop below below the threshold $m$ at a much later time $t\gg 20$. 
In order to reach such time-scales, a much higher resolution in both $x$ and $t$ will be needed, which is beyond the scope of the current paper.
\end{remark}

\begin{figure}[htb!]
 \includegraphics[width=0.49\textwidth]{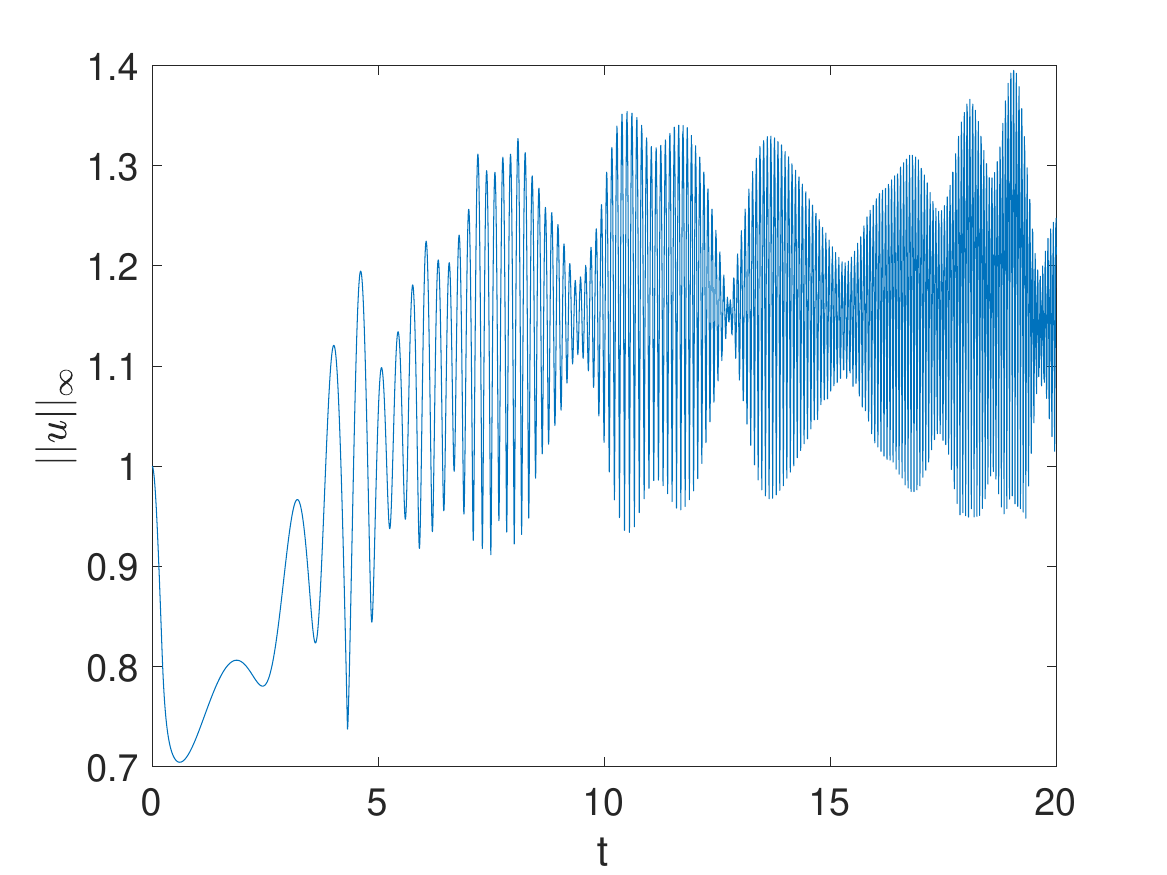}
 \includegraphics[width=0.49\textwidth]{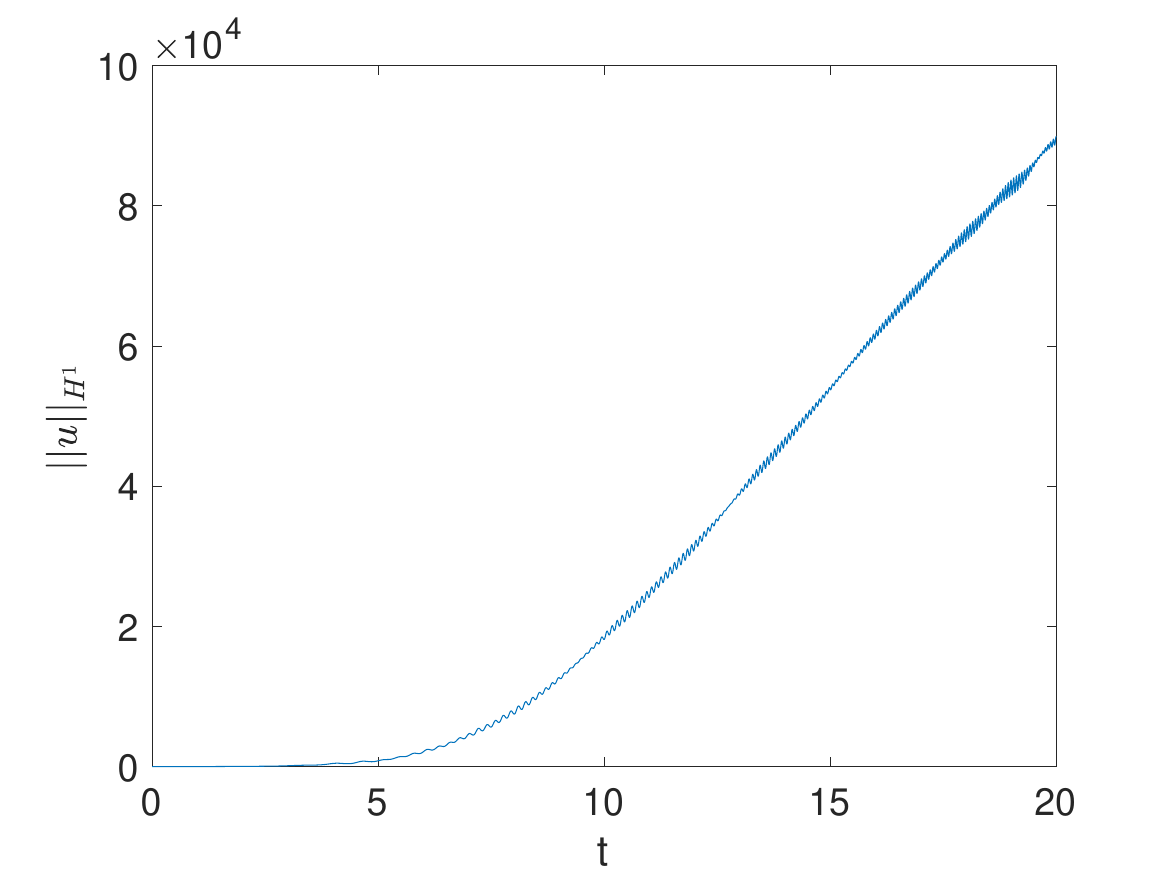}
 \caption{The $L^{\infty}$-norm and the $\dot {H^1}$-norm of the solution to the cubic fNLS \eqref{critfNLS} with $s=\frac{1}{4}$, $\eps=0.1$, and initial data \eqref{ininonlin}(a).}
 \label{fNLSs025t20inf}
\end{figure}

Further qualitative properties of the solution are shown in the next set of figures: In Fig.~\ref{fNLSs025t20inf} on the left, we show the $L^{\infty}$-norm of $u(t,\cdot)$ as a function of time. 
It is seen to be highly oscillatory, which is due to the oscillatory 
behavior of the solution and the fact that we numerically determine the values of $\|u(t, \cdot)\|_{L^\infty}$ only on grid 
points (the actual maximum of the solution might not be located on a 
point of the numerical grid). Extrapolating from the time interval shown, we expect the  
$L^{\infty}$-norm to remain bounded in time. 
The $\dot{H}^{1}$-norm of $u$, shown on the right of Fig.~\ref{fNLSs025t20inf}, is clearly growing in $t$ but seems to have an inflection point around $t\approx 20$. Thus 
we do not expect the appearance of a finite time singularity in this case.  

\begin{figure}[htb!]
 \includegraphics[width=0.49\textwidth]{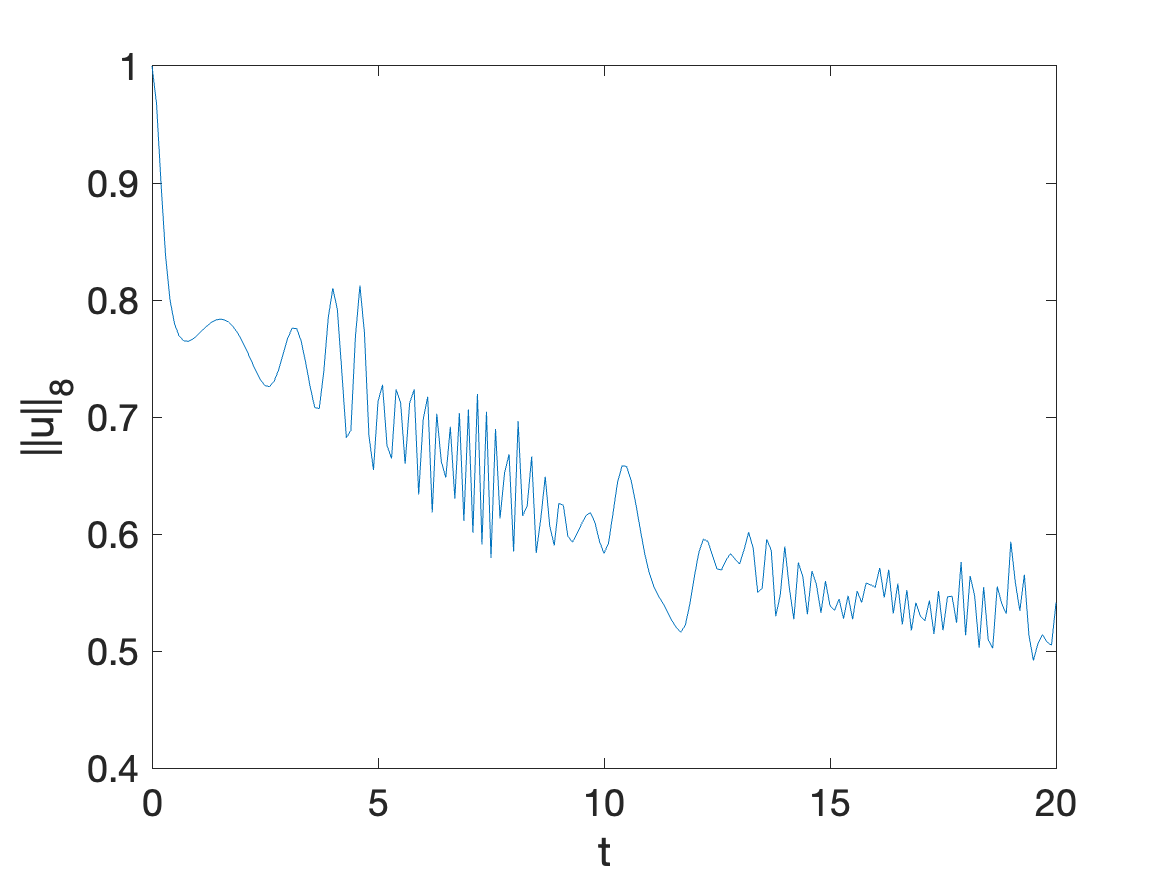}
 \includegraphics[width=0.49\textwidth]{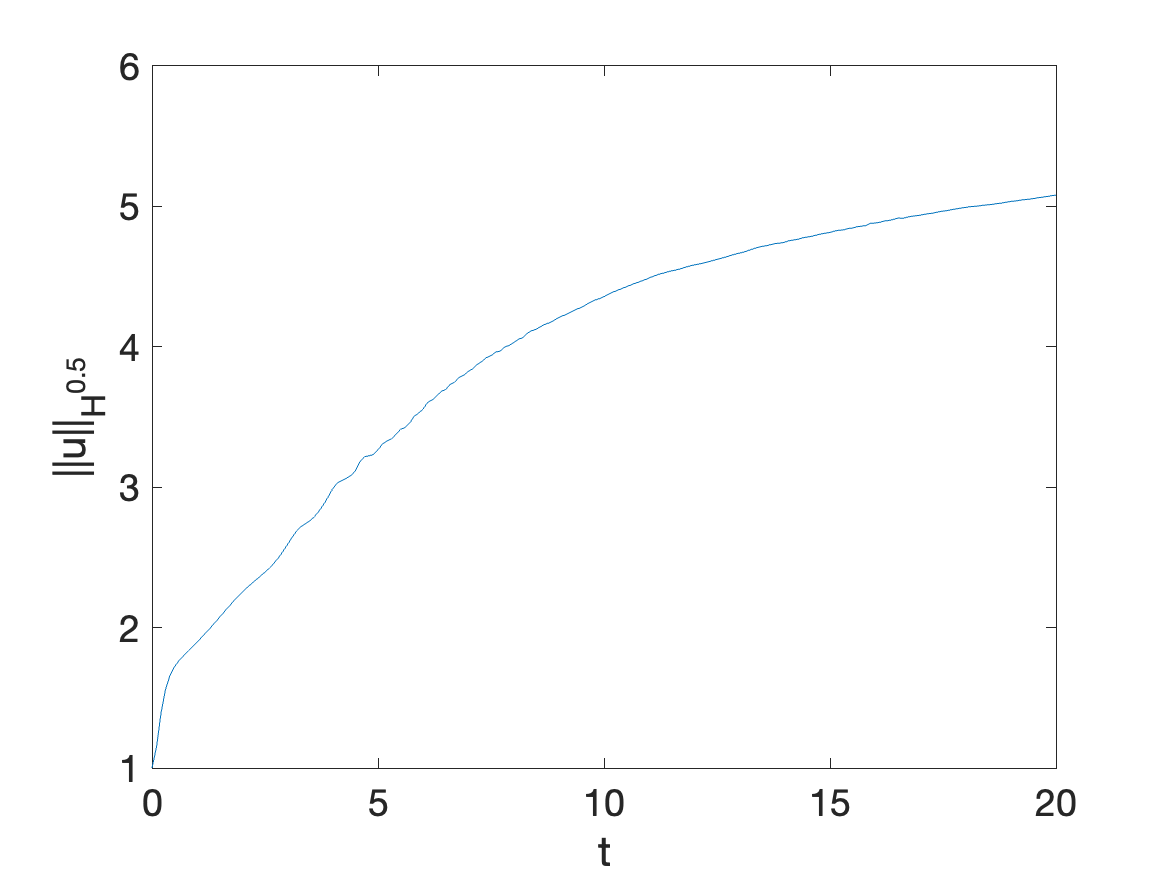}
 \caption{The $L^{8}$-norm and the $\dot {H}^{1/2}$-norm of the 
 solution to the cubic fNLS \eqref{critfNLS} with $s=\frac{1}{4}$, 
 $\eps=0.1$, and initial data \eqref{ininonlin}(a), both norms 
 normalized to 1 for $t=0$.}
 \label{fNLSs025t20L8}
\end{figure}

This belief is strengthened by the time-decay of higher order $L^{q}$-norms of $u(t, \cdot)$ with $q>2$. Shown on the left of 
Fig.~\ref{fNLSs025t20L8} is the case with $q=8$. In addition, we also plot the Sobolev-critical $\dot{H}^{1/2}$-norm of $u$ on the right 
of the same figure. After initially growing in $t$, it appears to become saturated and thus bounded. 

\begin{figure}[htb!]
 \includegraphics[width=0.49\textwidth]{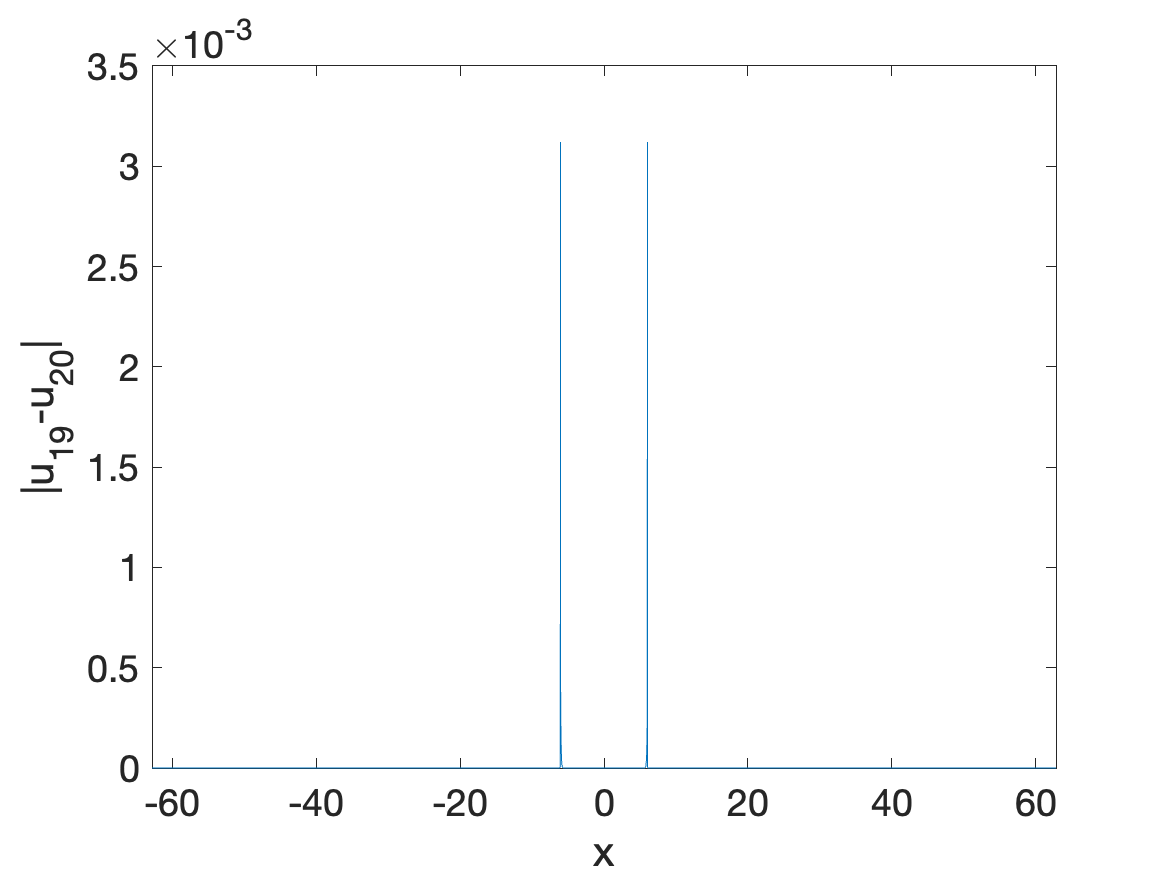}
 \includegraphics[width=0.49\textwidth]{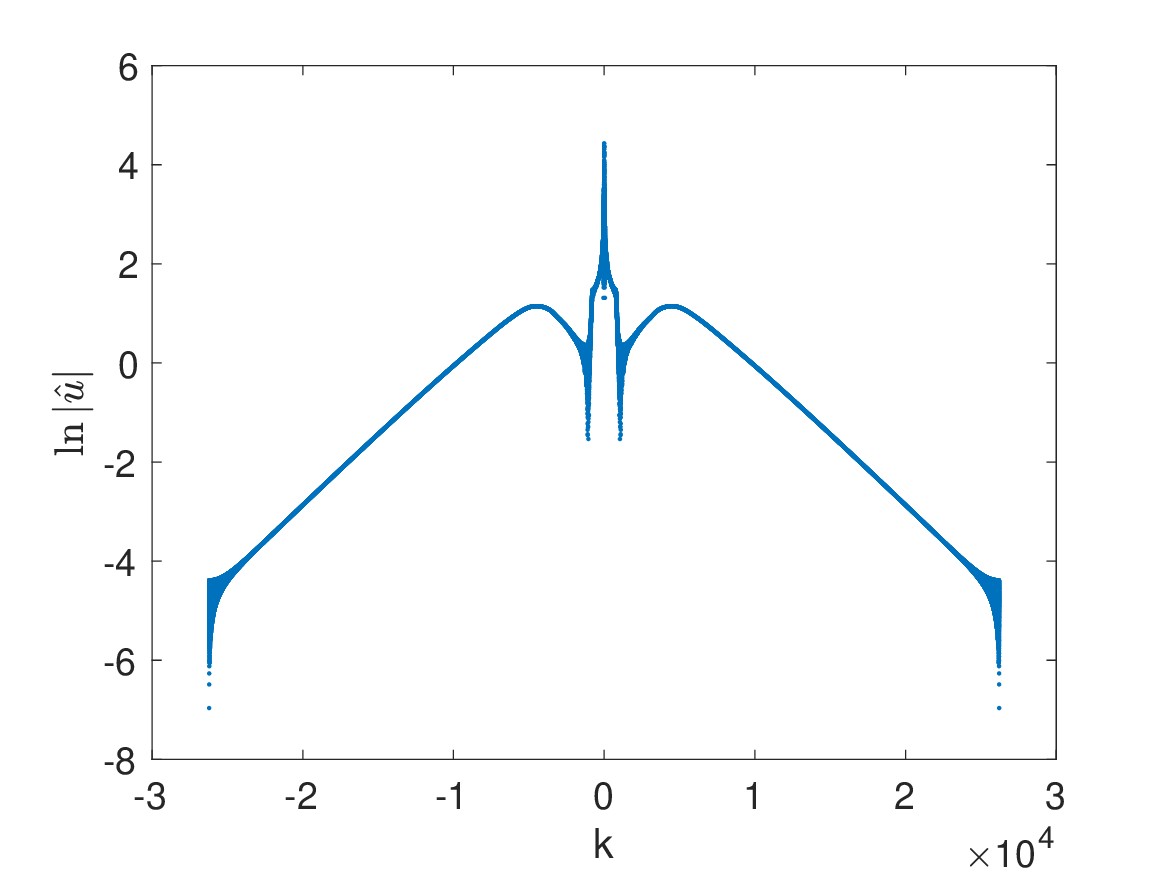}
 \caption{Left: difference between the numerical solution in 
 Fig.~\ref{fNLSs025t20} for $N=2^{19}$ and $N=2^{20}$ FFT modes.  
 Right: FFT coefficients of the solution from Fig.~\ref{fNLSs025t20} with $N=2^{20}$. }
 \label{fNLSs025t20fourierN20}
\end{figure}

To ensure that these findings are indeed reliable we redo our simulation: using varying time-steps, 
we first find that for our RK4 scheme the $L^{\infty}$-norm of the difference 
between the numerical solution obtained with $N_{t}=5\cdot10^{5}$ and the 
one with $N_{t}=10^{5}$ is of the order $2.4\cdot10^{-4}$, whereas the 
difference to the case with $N_{t}=5\cdot10^{4}$ is of the order of 
$6.7\cdot10^{-3}$. Thus for fixed spatial resolution, the numerical error 
due to the time discretization decreases as expected with $N_{t}$. 
In a similar way we can study the effect of spatial 
resolution. For fixed $N_{t}=10^{5}$, the difference between the 
numerical solution with $N=2^{19}$ and $N=2^{20}$ Fourier modes can be seen on 
the left of Fig.~\ref{fNLSs025t20fourierN20}. The FFT coefficients 
for the latter case at the final time are given on the right of the 
same figure. Visibly they decrease exponentially, confirming the 
regularity of the solution. In the appendix, we also compare our results with those obtained by two other numerical schemes to 
find that they agree within an $L^\infty$-difference of order $10^{-4}$.

\subsection{Further examples in the critical regime} 
\begin{figure}[htb!]
 \includegraphics[width=0.49\textwidth]{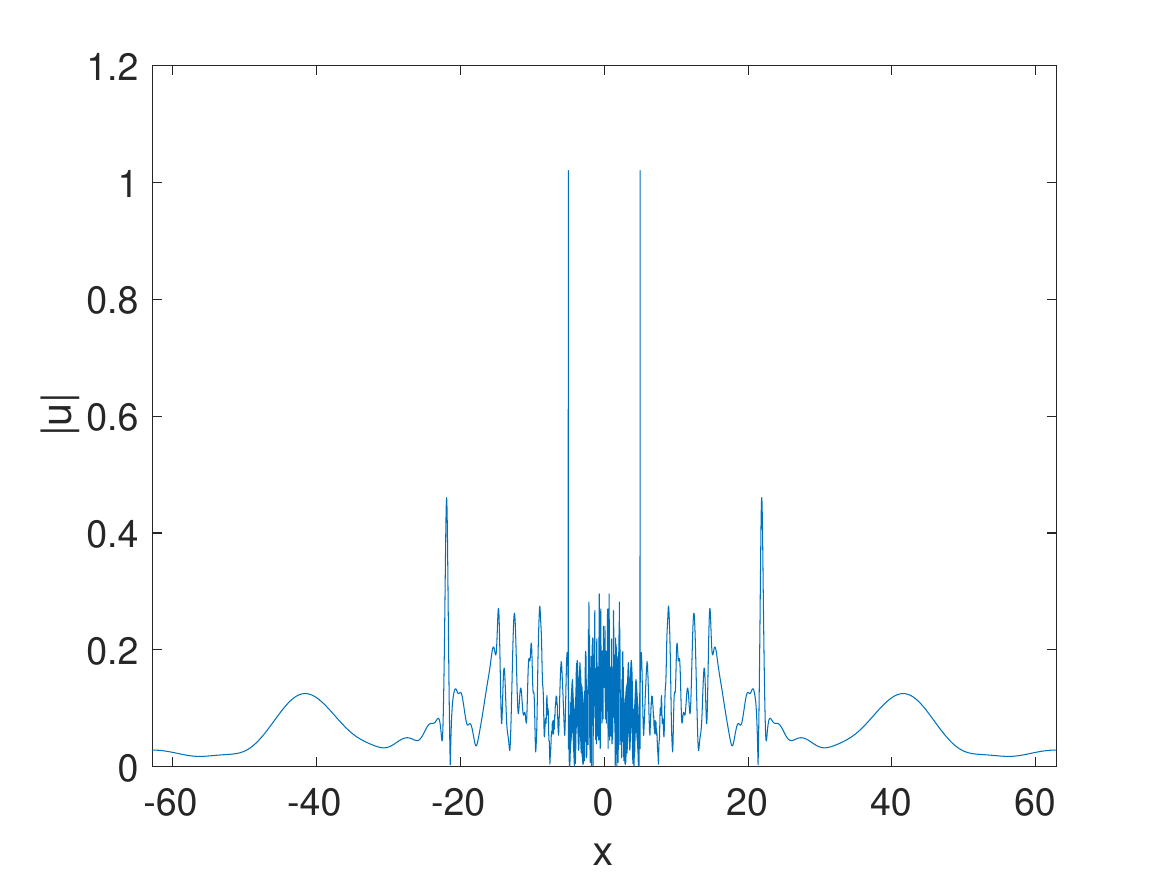}
 \includegraphics[width=0.49\textwidth]{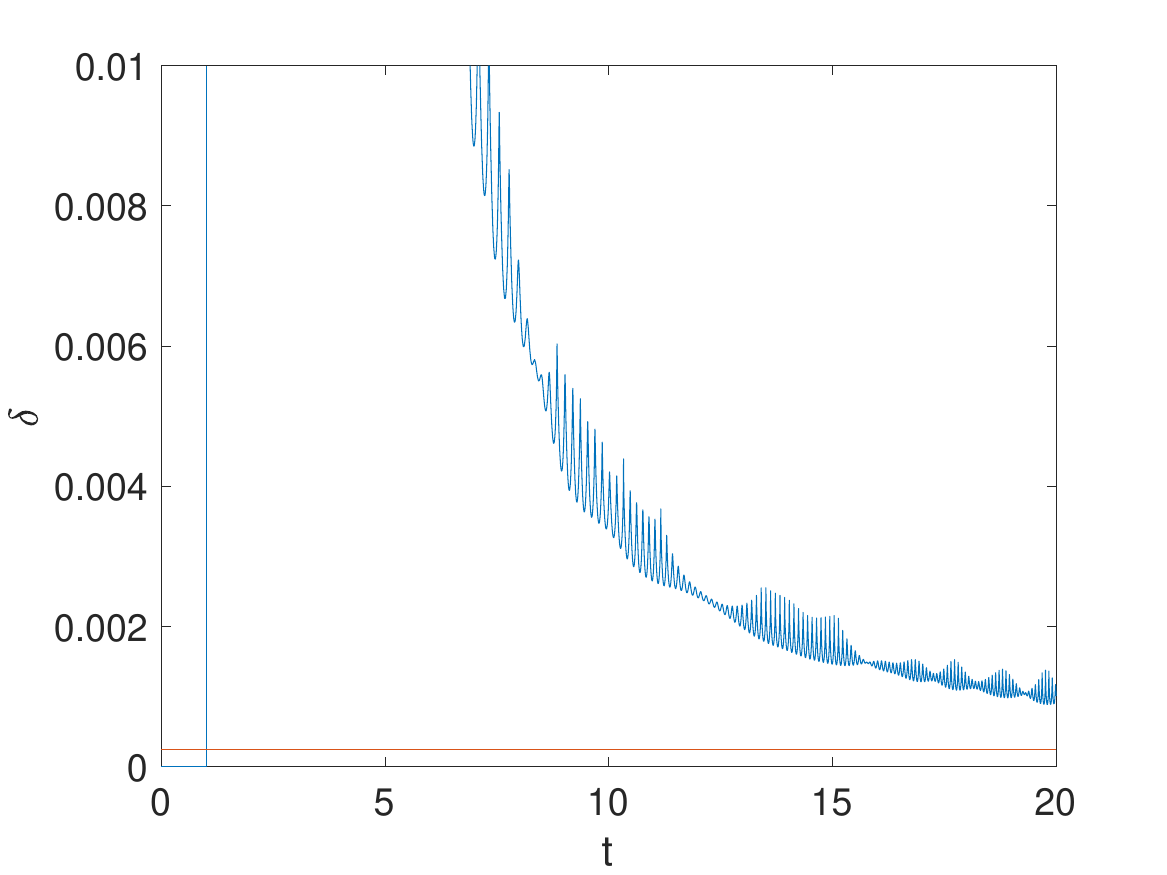}
 \caption{Left: solution to the cubic fNLS \eqref{critfNLS} with $s=\frac{1}{4}$, $\eps=0.1$, and initial data \eqref{ininonlin}(b) at $t=20$. Right: the 
 fitting parameter $\delta$ as a function of time.}
 \label{fNLSs025t20gauss}
\end{figure}

The situation remains qualitatively similar in the case of Gaussian initial data \eqref{ininonlin}(b). 
The corresponding modulus of the solution $u$ at time $t=20$ is shown 
on the left of Fig.~\ref{fNLSs025t20gauss}. 
Notably the bulk of $|u|$ is more localized than in the previous case. The function $\delta(t)$ 
is shown on the right hand side of the same figure. As one can see, $\delta(t)$ decays slower than in the previous case and remains clearly above the threshold $m$.  

\begin{figure}[htb!]
 \includegraphics[width=0.49\textwidth]{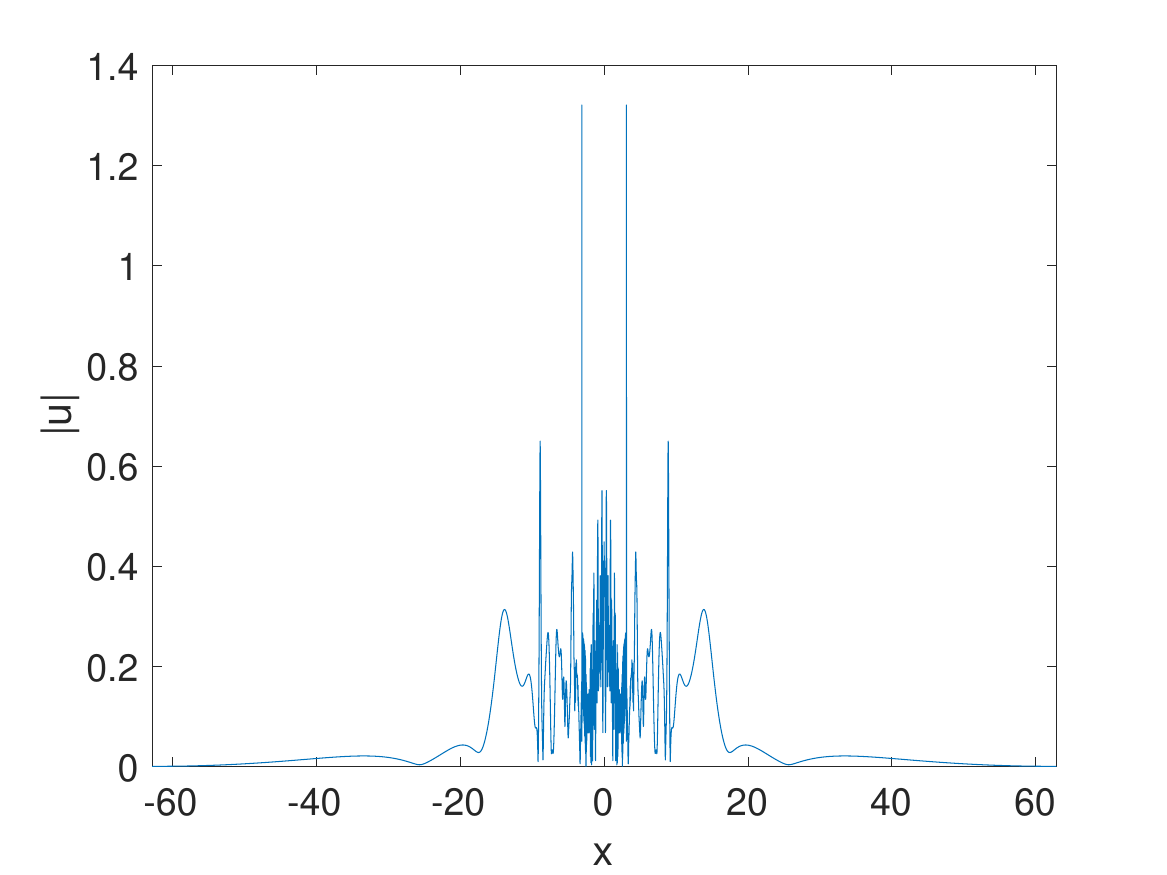}
 \includegraphics[width=0.49\textwidth]{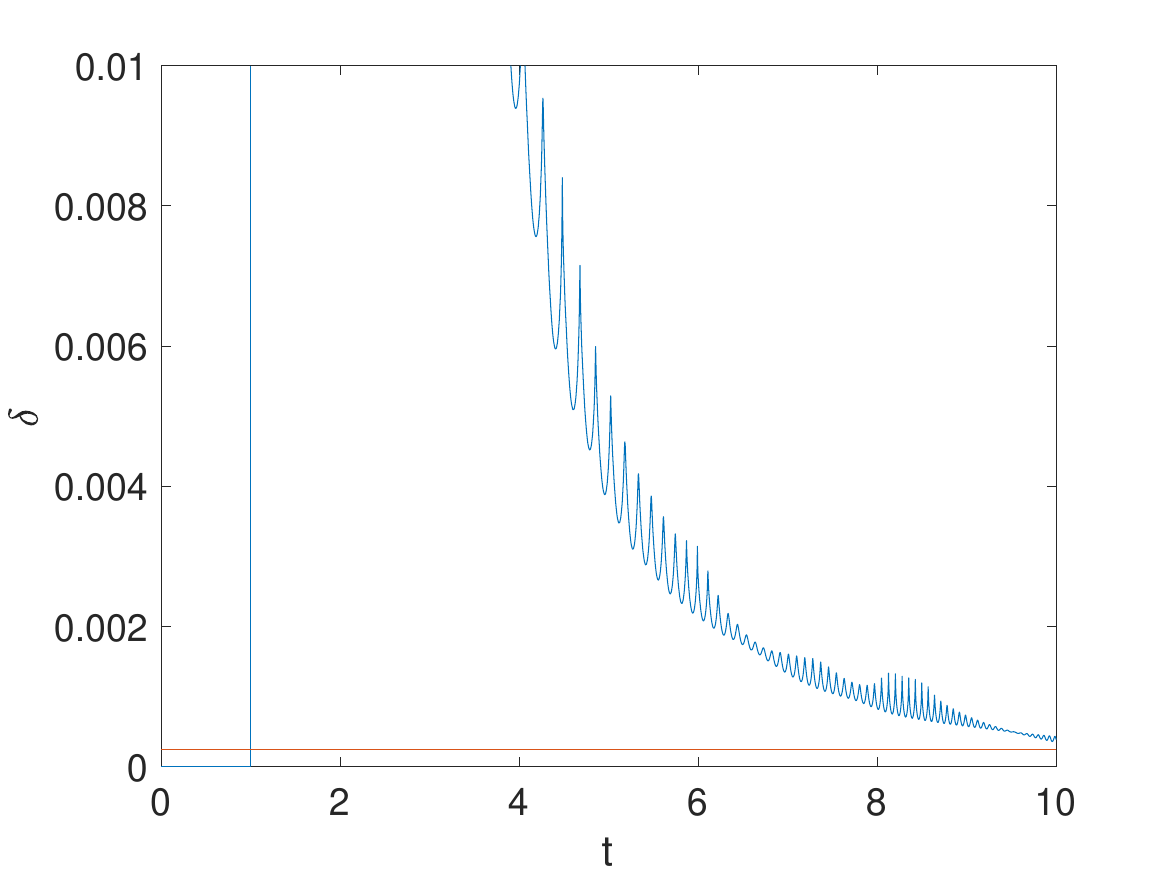}
 \caption{Left: solution to the cubic fNLS \eqref{critfNLS} with $s=\frac{1}{4}$, $\eps=0.1$, and initial data \eqref{ininonlin}(c) at $t=10$. Right: the 
 fitting parameter $\delta$ as a function of time.}
 \label{fNLSs025t10ex4}
\end{figure}

For the even faster decaying (super Gaussian) initial data \eqref{ininonlin}(c), the solution at the time 
$t=10$ is shown on the left of Fig.~\ref{fNLSs025t10ex4}. 
There are 
again rapid oscillations located between two clearly visible sharp peaks, but the solution 
is even more localized than before. In addition, the parameter $\delta(t)$ decreases more rapidly than 
in the previous two cases.


\section{Solution of the fractional NLS: the supercritical regime}\label{sec:supcrit}

\subsection{Cubic nonlinearity} We now turn to the supercritical regime, and first consider the 
following fNLS with cubic nonlinearity:
\begin{equation}\label{supfNLS}
i \eps \partial_t u  = (- \eps^2 \Delta )^{1/5} u + |u|^{2} u, \quad u_{\vert{t=0}} = \upsilon.
\end{equation}
Recall that for cubic nonlinearities the energy critical power is $s_\ast=\frac{1}{4}$, hence choosing $s=\frac{1}{5}$ puts us firmly in an energy supercritical regime. 
As initial data we shall again choose one of the functions given in \eqref{ininonlin}.

First, in the case where $\upsilon=\mbox{sech}(x)$, the resulting (weakly turbulent) time-evolution of $|u|$ is shown in the next figure. 
\begin{figure}[htb!]
 \includegraphics[width=0.8\textwidth]{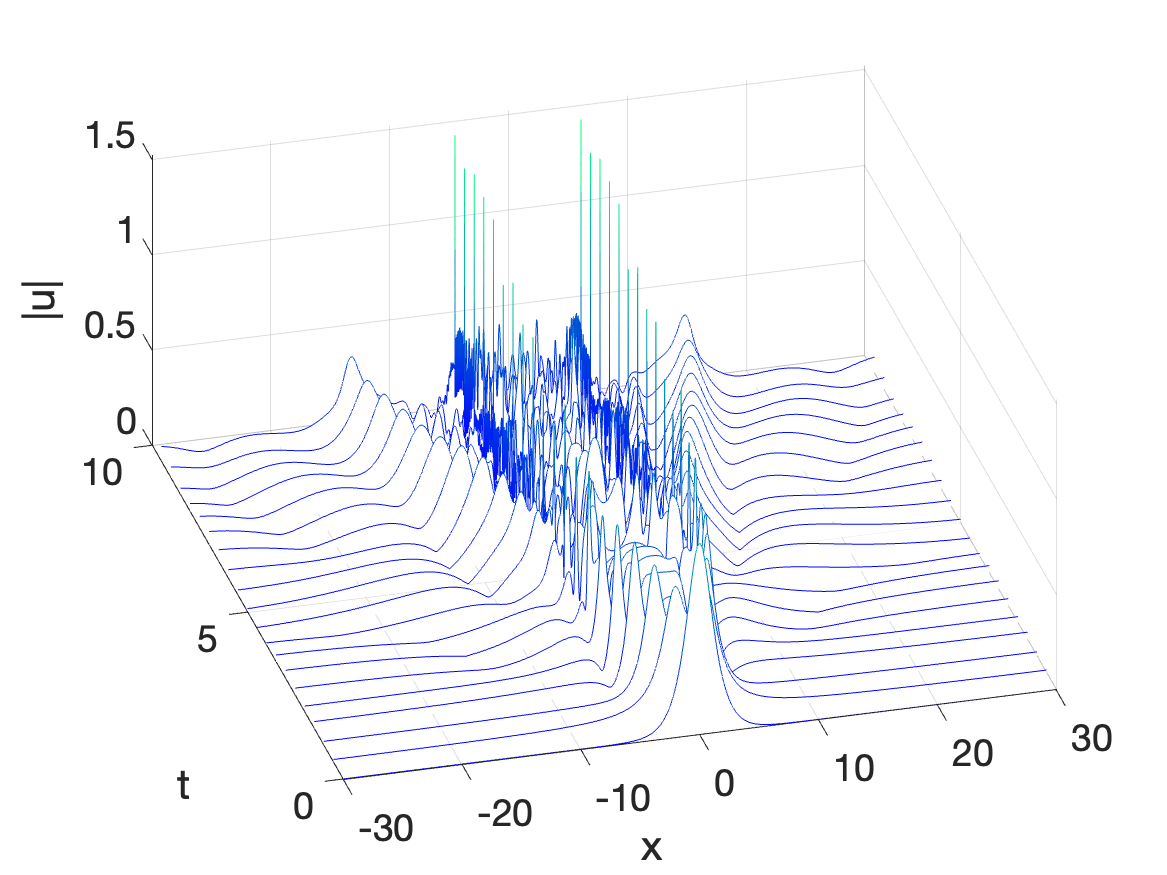}
\caption{Solution to the cubic fNLS 
\eqref{supfNLS} with $\eps=0.1$, $s=\frac{1}{5}$ and initial data \eqref{ininonlin}(a). }
\label{fNLSs02t98water}
\end{figure}

\begin{figure}[htb!]
 \includegraphics[width=0.49\textwidth]{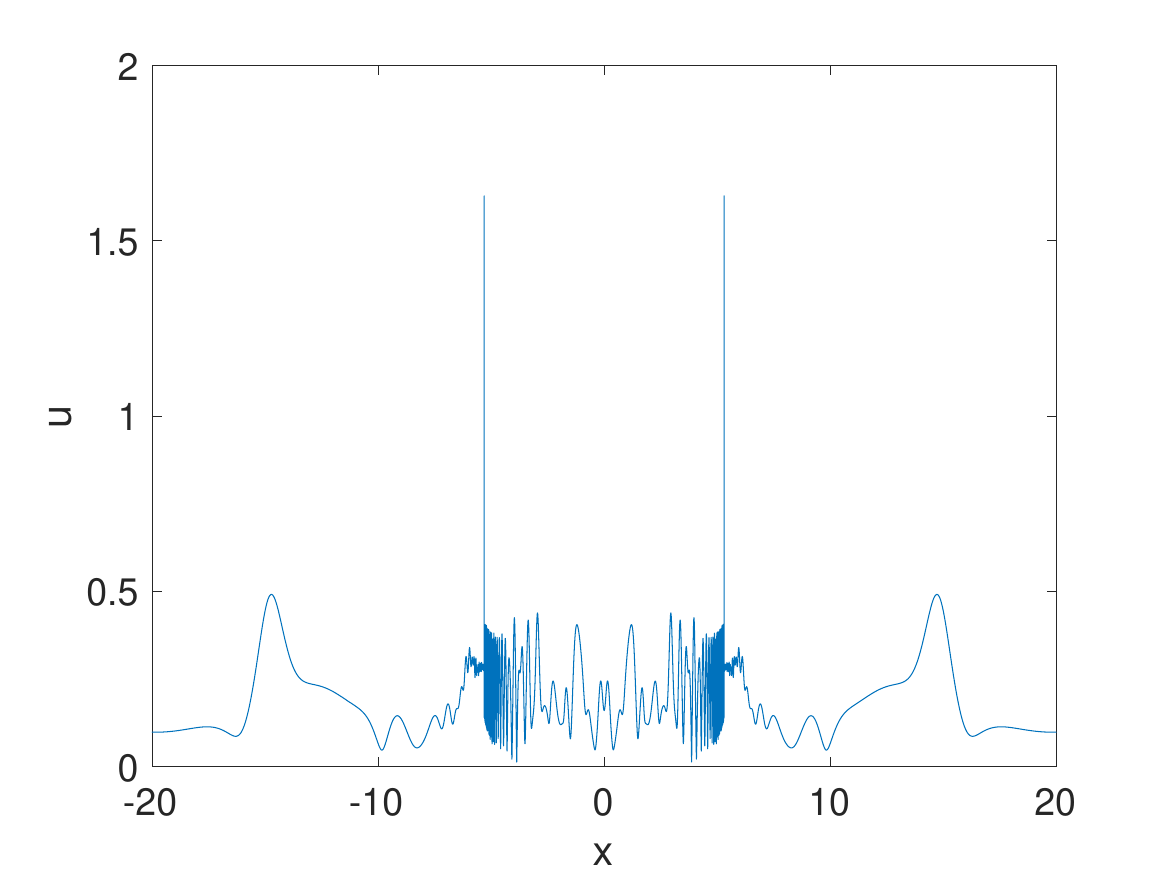}
 \includegraphics[width=0.49\textwidth]{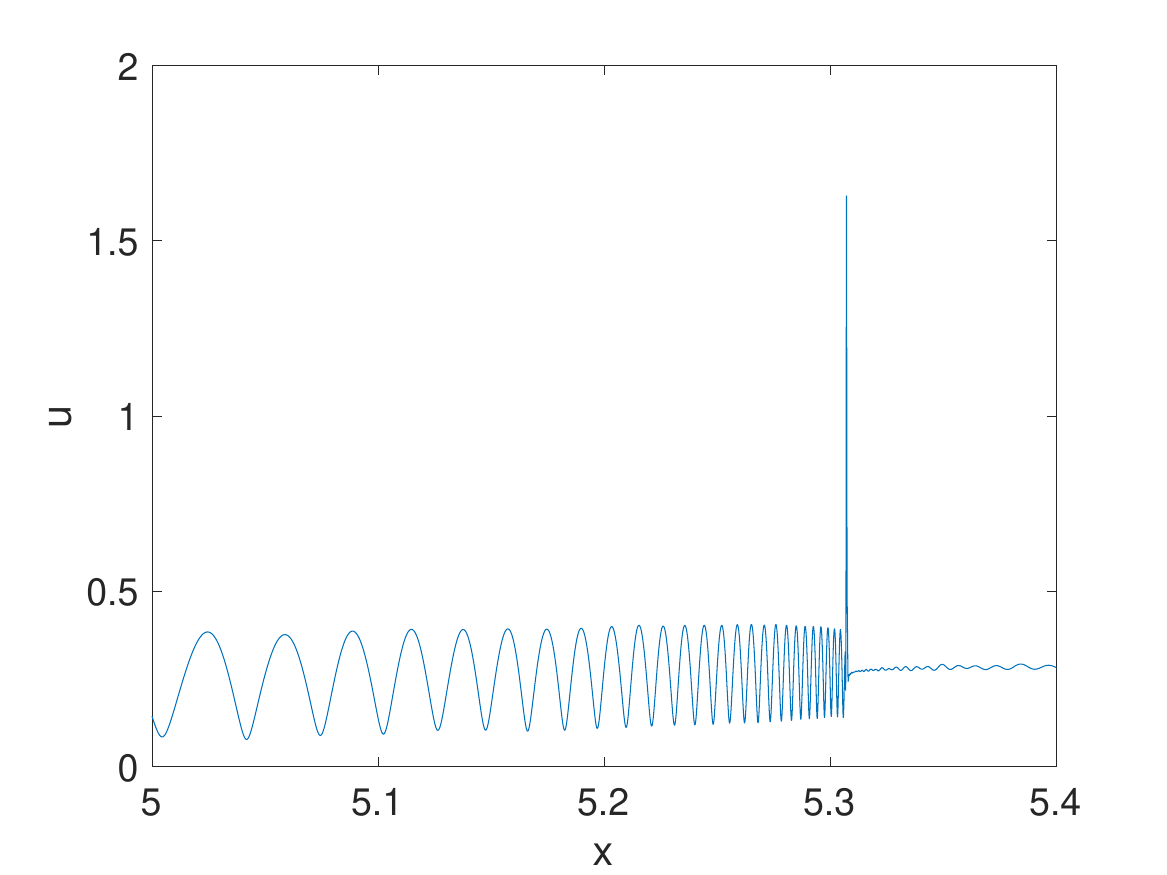}
 \caption{Left: solution to the cubic fNLS \eqref{supfNLS} with $\eps=0.1$, $s=\frac{1}{5}$ and initial data \eqref{ininonlin}(a) at the final time $t_{\rm f}=9.8$. Right: a close-up 
 near one of the peaks. }
 \label{fNLSs02t98}
\end{figure}

In this case, the numerical value of $\delta(t)$ vanishes at time $t_{\rm f}= 9.8$. 
The code is stopped at this point, and the solution at the final time can be seen in Fig.~\ref{fNLSs02t98}. 

As in the energy critical case there appear to be 
strong peaks at the left and right front of the two humps  
with rapid oscillations in between. However, in contrast to the previous (sub-)critical situation, these 
oscillations become more and more ``compressed" close to the peaks. This is 
clearly visible on the right of the same figure, which 
shows a close-up of the singular region around $x\simeq 5.31$. 

The fitting of 
the FFT coefficients according to \eqref{fourasymp} yields a value 
of $\mu\simeq -\frac13$ (numerically, we obtain $\mu=-0.327\dots$). This indicates a blow-up of the $L^{\infty}$-norm of $u$, 
since, according to the theory presented in Section \ref{sec:singtr}, we expect
\[
|u(t_{\rm f}, x)| \sim |x-x_{0}|^{-\frac13}\ \text{where $x_0\simeq \pm 5.31$}.
\]
However, one might need to take this with a grain of salt (cf. the discussion at the end of Section \ref{sec:singtr}). 
A finite amplitude $|u(t_{\rm f}, x_0)|$ which exhibits blow-up within some higher order $L^q$-norm with $q> 2$, or some higher order $H^s$-norm with $s>\frac12$, 
seems equally likely, and hard to decide numerically. 
One should also note that the relatively small value of $|\mu|\simeq \frac{1}{3}$ implies that the singularity at $x_0$ would yield a (very) large $L^{\infty}$-norm only if $x_0$ were by chance {\it exactly} on one of 
the grid points.

\begin{figure}[htb!]
 \includegraphics[width=0.49\textwidth]{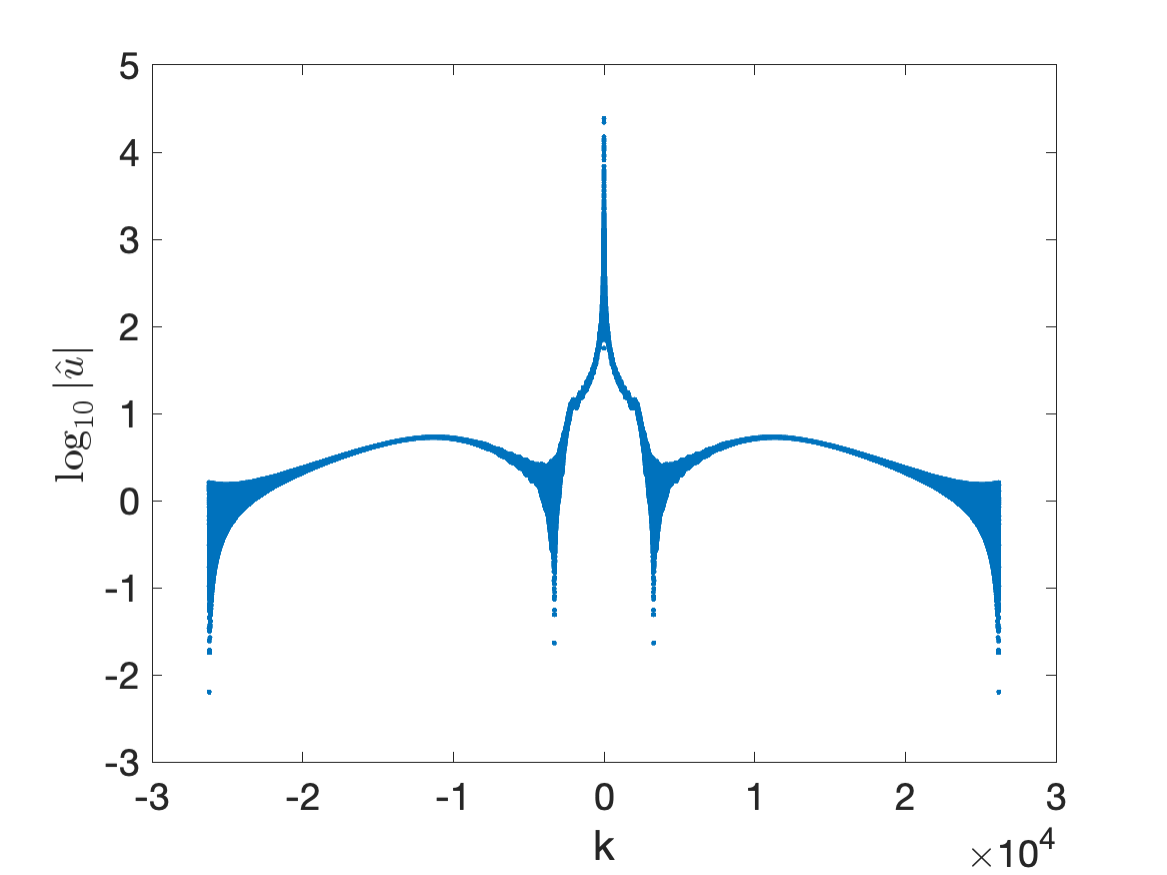}
 \includegraphics[width=0.49\textwidth]{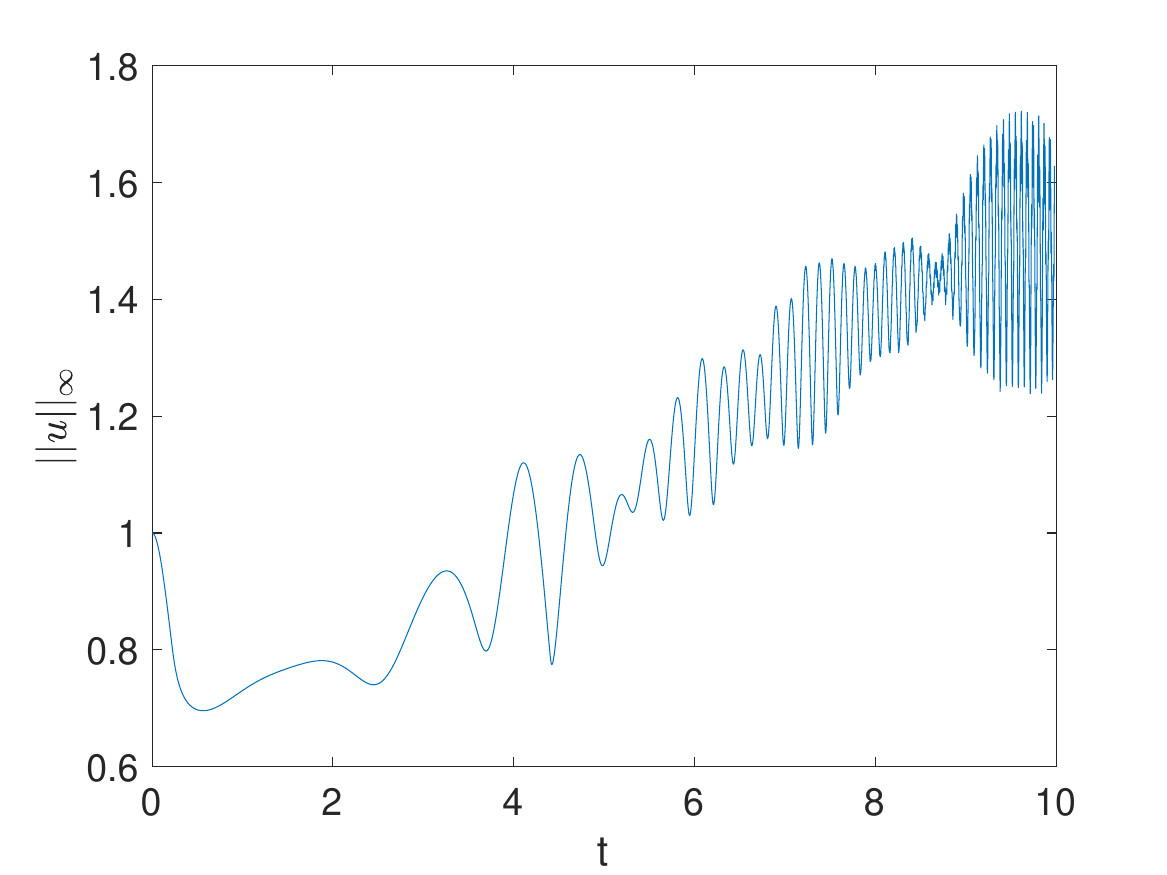}
 \caption{Left: FFT coefficients of the solution from Figure \ref{fNLSs02t98}. Right: the $L^{\infty}$-norm of $u$ as a function of time.  }
 \label{fNLSs02t98fourier}
\end{figure}

The FFT coefficients within $\widehat u$ at time $t_{\rm f}=9.8$ are shown on the left of 
Fig.~\ref{fNLSs02t98fourier}. On the right of the same figure we also show the
$L^{\infty}$-norm of the solution $u$ as a function of time. It can be seen that $\| u(t, \cdot) \|_{L^\infty}$ is slowly growing in time and again highly oscillatory. 

\begin{figure}[htb!]
 \includegraphics[width=0.49\textwidth]{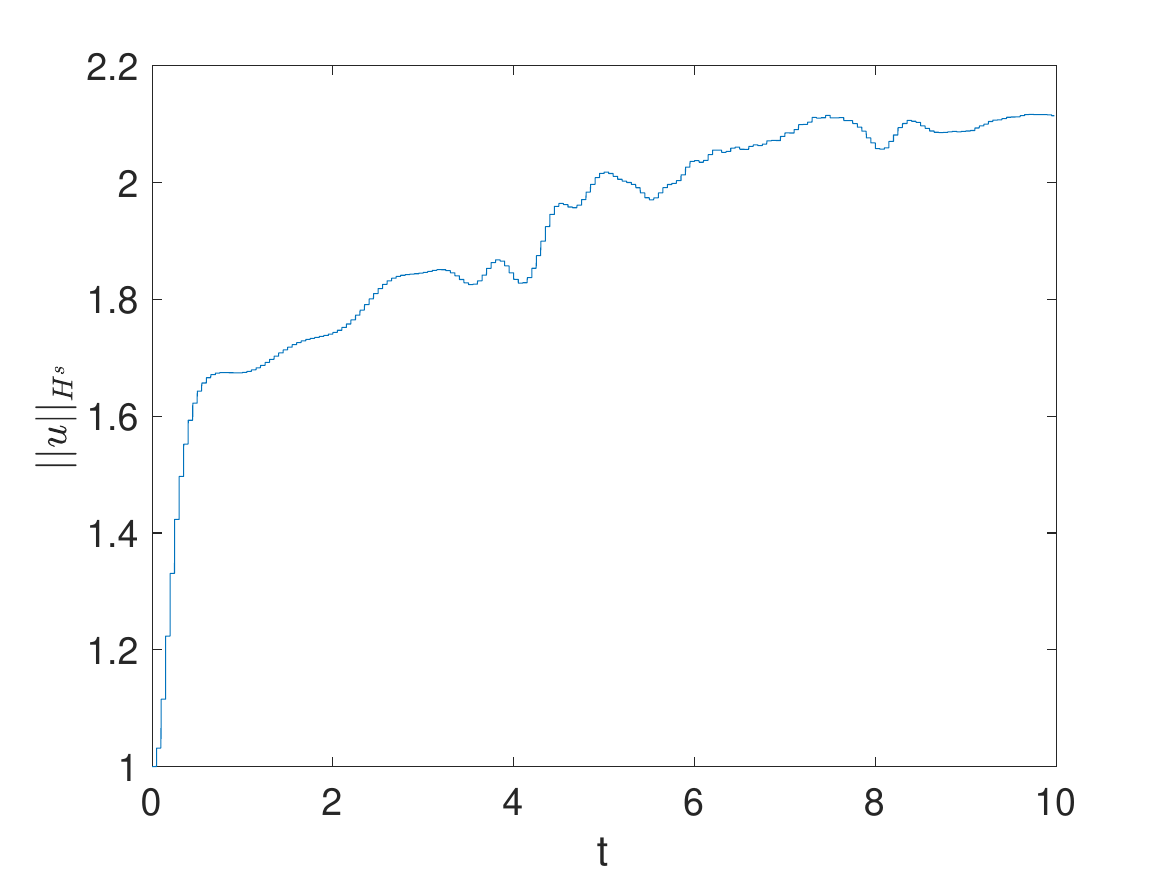}
 \includegraphics[width=0.49\textwidth]{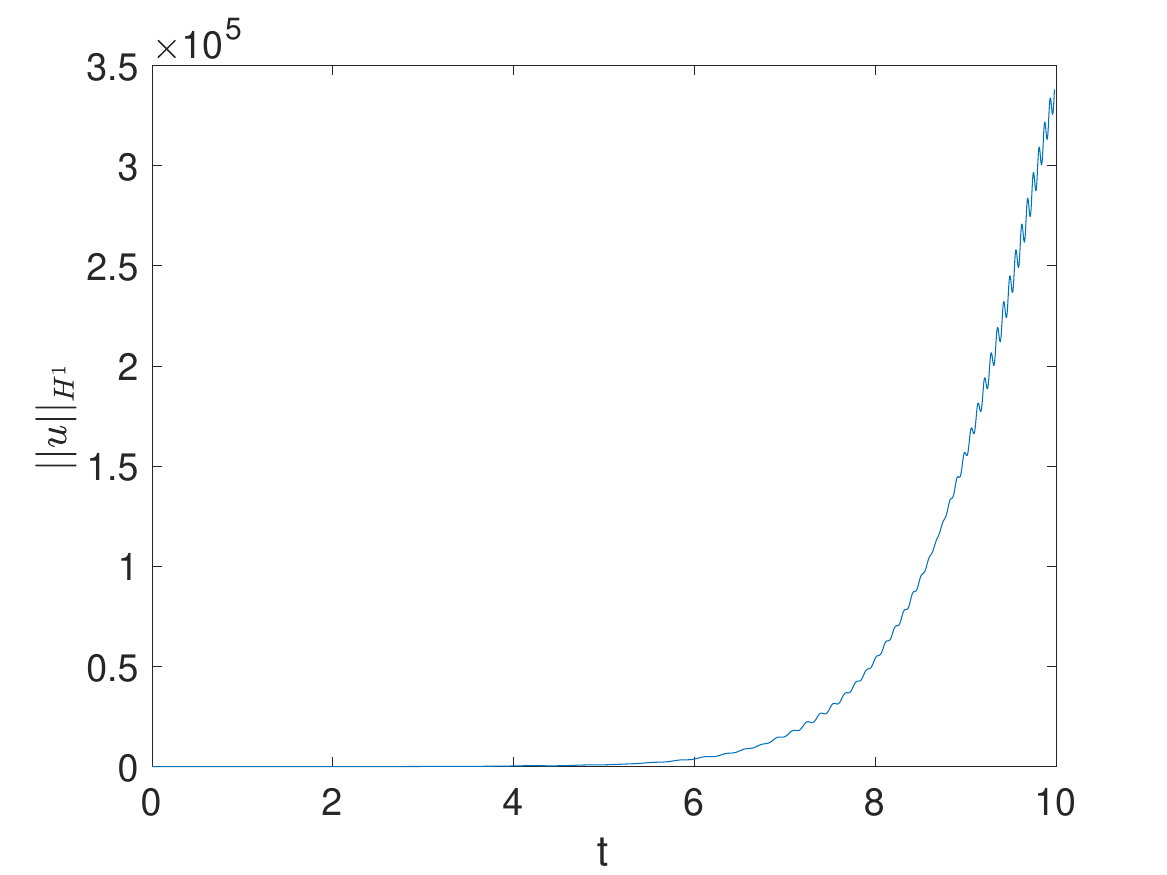}
 \caption{The $\dot{H}^{s}$-norm (left) and the $\dot{H}^1$-norm (right) of the solution to the cubic fNLS \eqref{supfNLS} with $\eps=0.1$, $s=\frac{1}{5}$ and initial data \eqref{ininonlin}(a). }
 \label{fNLSs02t98norm}
\end{figure}

In Fig.~\ref{fNLSs02t98norm} we show the behavior of the (energy-controlled) $\dot{H}^{s}$-norm of $u$ 
and of the supercritical $\dot{H}^{1}$-norm, both as functions of time. 

\begin{figure}[htb!]
 \includegraphics[width=0.49\textwidth]{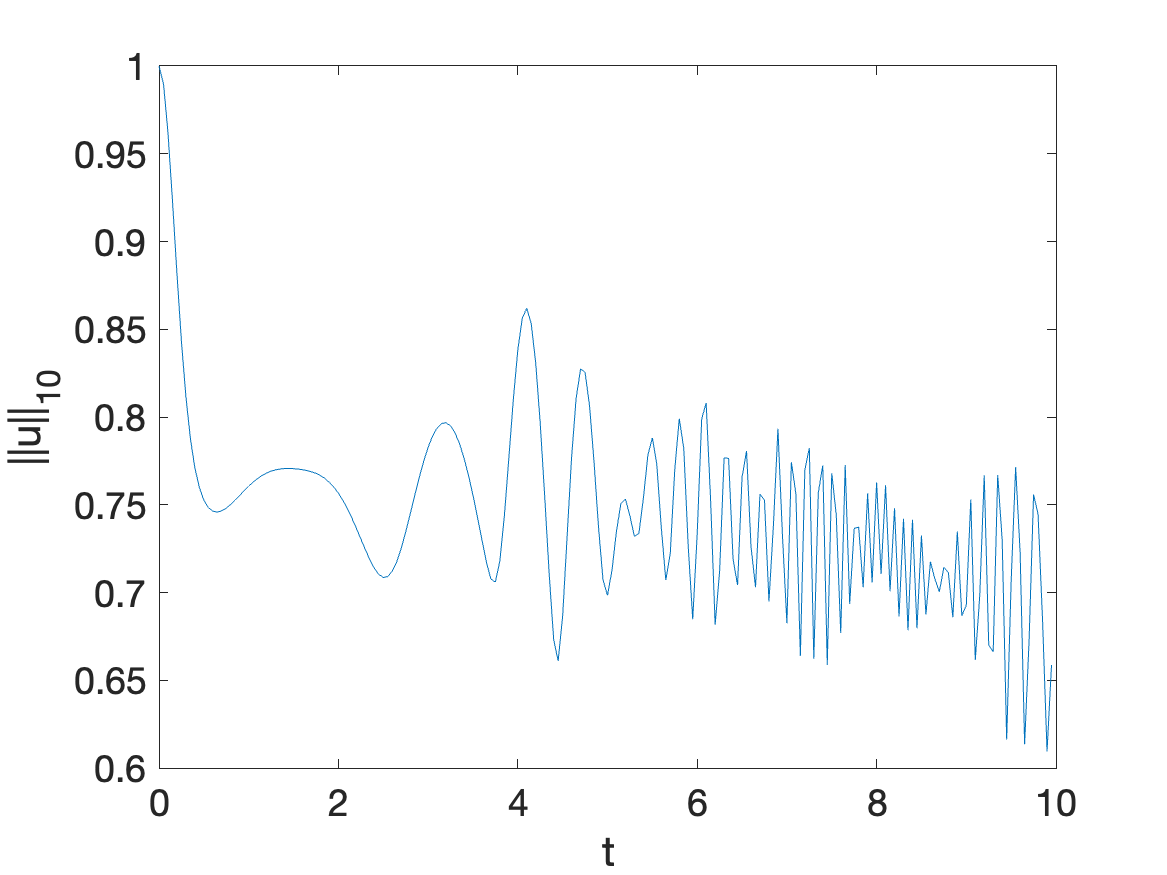}
 \includegraphics[width=0.49\textwidth]{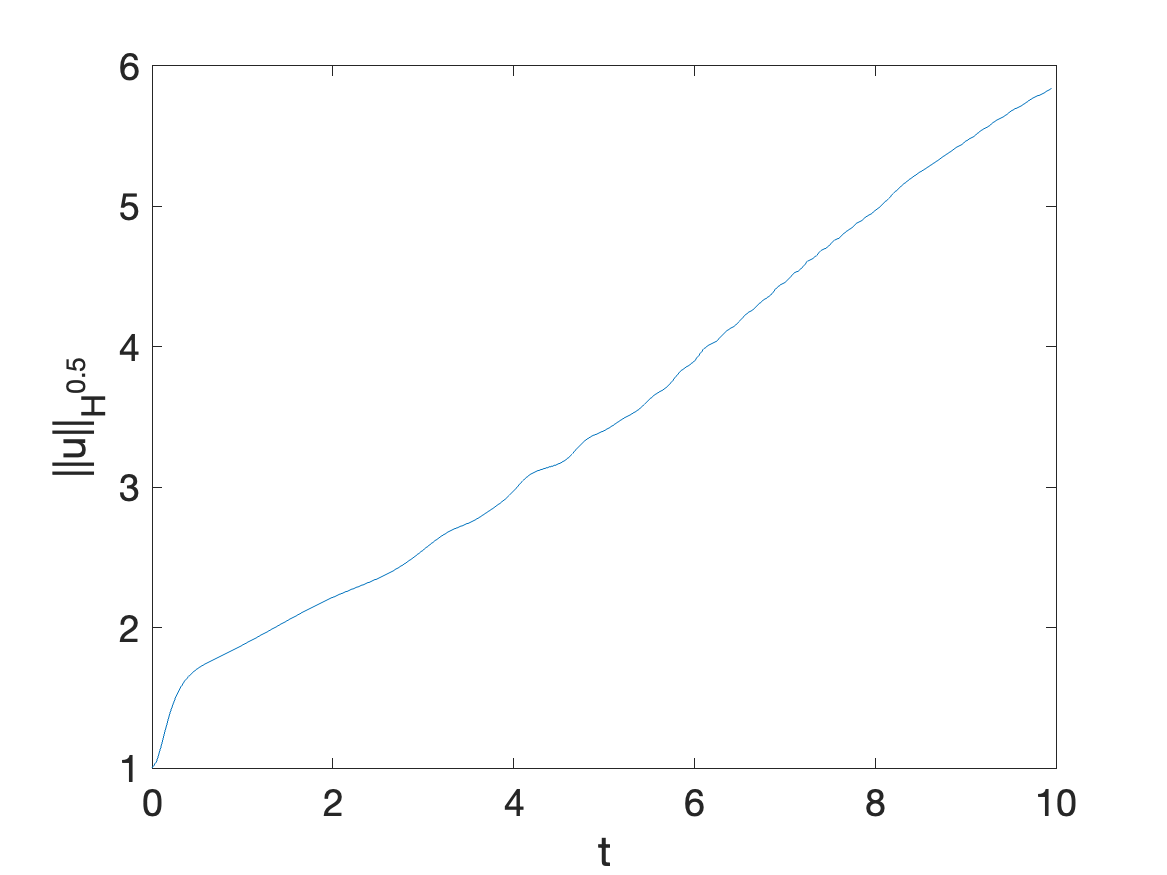}\\
 \caption{The $L^{10}$-norm (left) and the $\dot{H}^{1/2}$-norm 
 (right) of the solution to the cubic fNLS \eqref{supfNLS} with 
 $\eps=0.1$, $s=\frac{1}{5}$ and initial data \eqref{ininonlin}(a), 
 both normalized to 1 for $t=0$.
}
\label{fNLSs02t98L10}
\end{figure}

It can be seen 
that the former remains bounded, due to (numerically almost) energy 
conservation, whereas the latter is clearly rapidly growing, indicating a finite-time blow-up. Other supercritical Sobolev norms $H^s$ with $s\ge 1$ (not shown here) behave similarly.
In contrast, all $L^{q}$-norms of $u(t, \cdot)$ with $2\le q\leq 10$ remain bounded as functions of time, cf. Fig.~\ref{fNLSs02t98L10} which 
shows the case for $q=10$. On the right of the same 
 figure we show the Sobolev critical $\dot{H}^{1/2}$-norm, which is seen to be moderately growing. 
 
Overall, the numerical pictures suggest the {\it formation of a small-scale coherent structure} associated with the transfer of energy to arbitrarily small spatial regions: a “spike” of large amplitude, 
whose width seems to go to zero, approaching a point discontinuity at some finite time $T\gtrsim t_{\rm f}$.
All of this is combined with a region of high frequency oscillations within $u$, suggesting the appearance of a singularity in {\it both} the amplitude $|u|$ and the phase $\varphi = \arg u$.

\begin{remark} Very recently, numerical and analytical evidence for a similar phenomenon, called {\it -3/2-cascade}, has been found in the 
case of fully resonant Hamiltonian systems with cubic nonlinearity, cf. \cite{Bi, BiGe}. The latter can be seen as a non-dispersive toy model, obtained 
from the fNLS by ignoring the fractional Laplacian and retaining only cubically resonant 
frequencies $k\in \Z$, which satisfy the relation: $|k_1|^{2s}-|k_2|^{2s}+|k_3|^{2s} = |k|^{2s}$ for $k_1-k_2+k_3 = k$. 
In particular, the figures in \cite{Bi, BiGe} show the formation of a ``spike", which very much resembles our findings, minus the additional high frequency oscillations 
(which presumably are due to non-resonant frequencies).
\end{remark}

\subsection{Additional simulations in the supercritical cubic case} 
\begin{figure}[htb!]
 \includegraphics[width=0.49\textwidth]{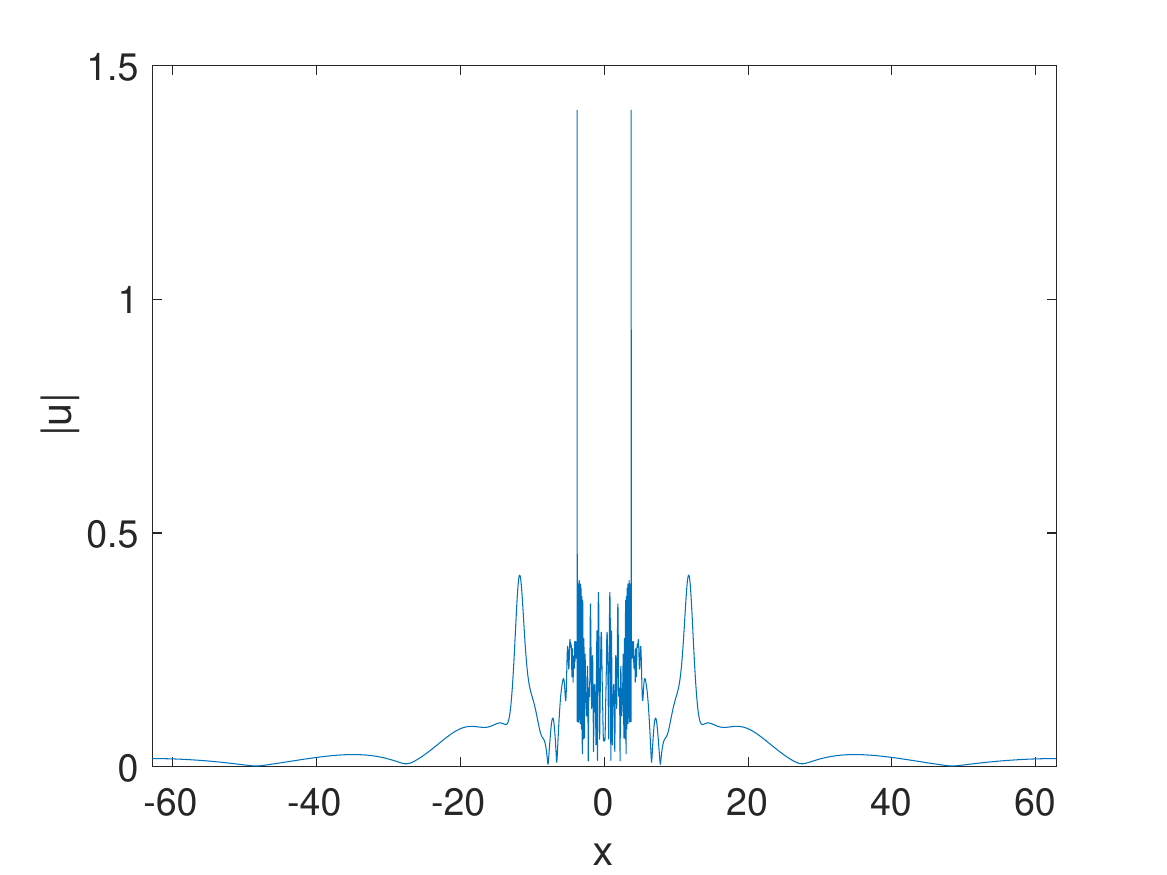}
 \includegraphics[width=0.49\textwidth]{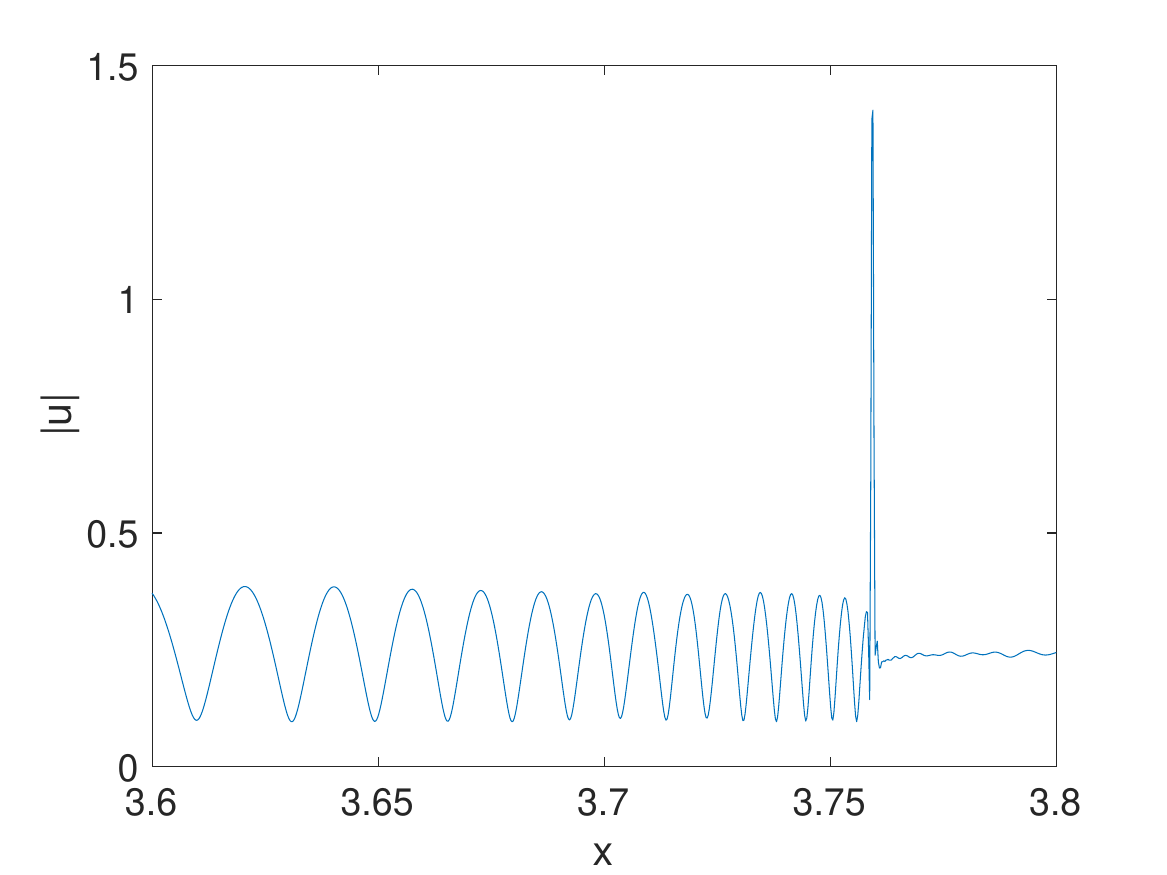}
 \caption{Left: solution to the cubic fNLS \eqref{supfNLS} with $\eps=0.1$, $s=\frac{1}{5}$ and initial data \eqref{ininonlin}(b) at the final time $t_{\rm f}=9.59$. 
 Right: a close-up 
 near one of the peaks. }
 \label{fNLSs02t959}
\end{figure}

\begin{figure}[htb!]
 \includegraphics[width=0.49\textwidth]{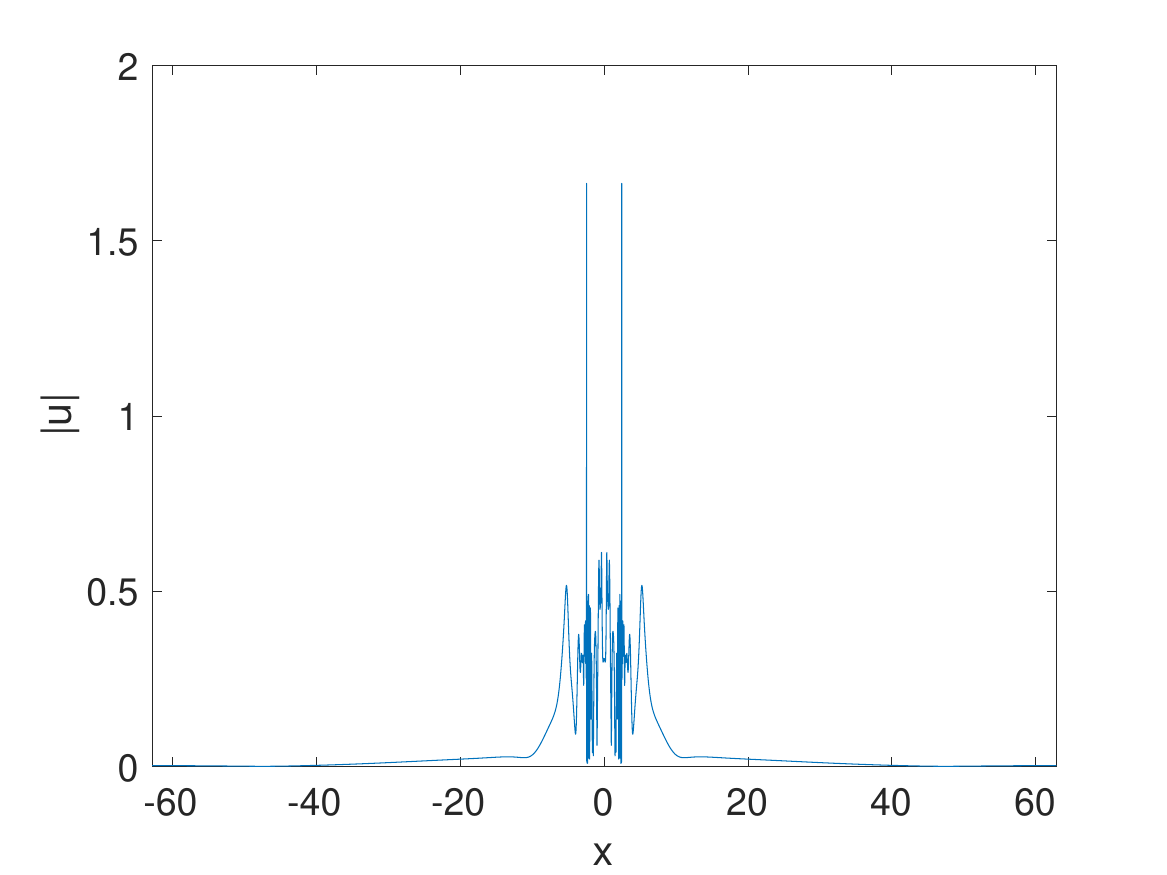}
 \includegraphics[width=0.49\textwidth]{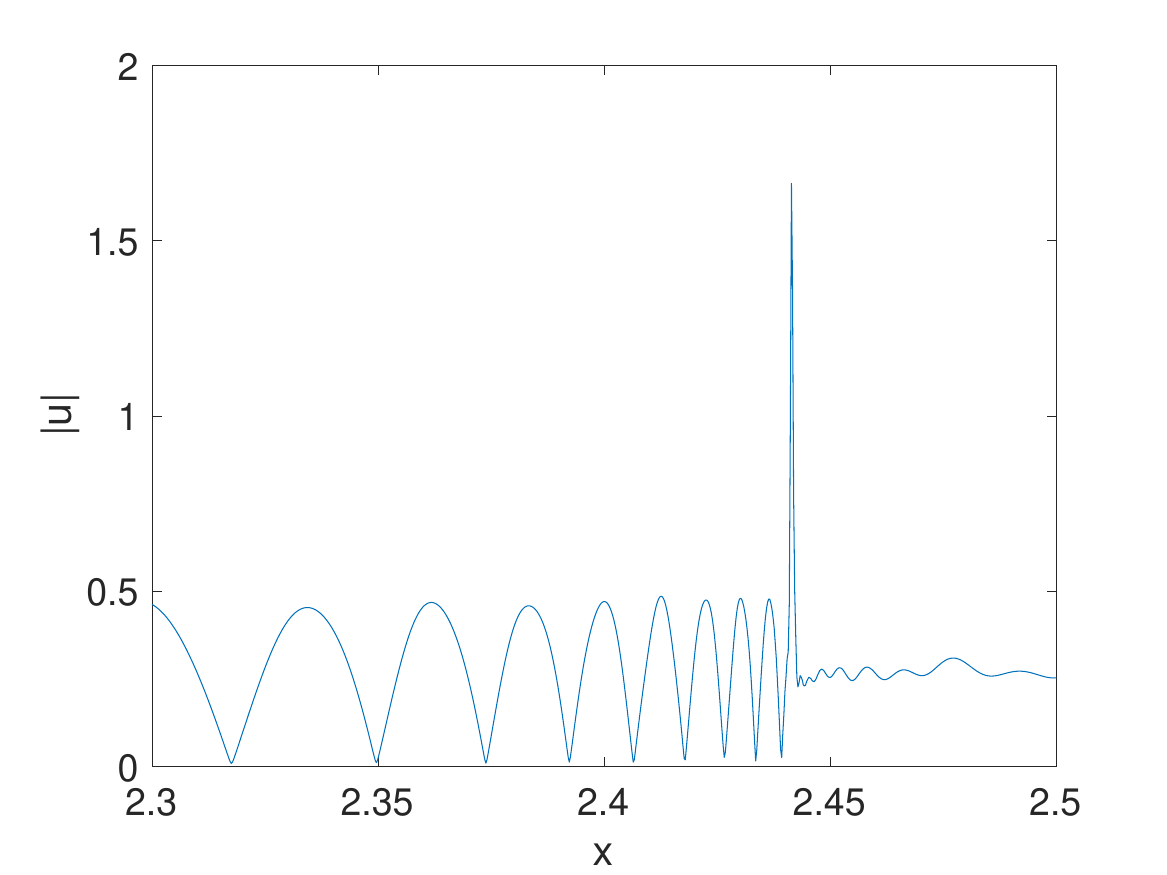}
 \caption{Left: solution to the cubic fNLS \eqref{supfNLS} with $\eps=0.1$, $s=\frac{1}{5}$ and initial data \eqref{ininonlin}(c) at the final time $t_{\rm f}=9.59$. 
 Right: a close-up 
 near one of the peaks. }
 \label{fNLSs02t546}
\end{figure}

For Gaussian initial data \eqref{ininonlin}(b), the code is stopped at time $t_{\rm f}= 9.59$,  as $\delta =0 $ at this time. We thereby obtain a 
fitted value of $\mu\simeq -0.33\dots$. The solution at the final recorded time is 
shown in Fig.~\ref{fNLSs02t959}. When compared to Fig.~\ref{fNLSs02t98}, the solution is qualitatively similar but notably more confined. A 
close-up of the solution near one of the peaks on the right of 
Fig.~\ref{fNLSs02t959} again shows compression of the 
oscillations near the presumed singularity at $x_0\simeq 3.76$. 
Finally, for super-Gaussian initial data \eqref{ininonlin}(c), we find that numerically $\delta = 0$ at time $t_{\rm f}= 5.46$. The fitted value for $\mu$ is again compatible with $\mu \simeq -\frac13$ (we obtain 
$\mu = -0.32\dots$). The solution at the final time is shown in Fig.~\ref{fNLSs02t546} on 
the left, with a close up near the presumed singularity on the right. 

The above runs have been repeated with $N=2^{20}$ FFT modes, and all of our results agree within the estimated numerical accuracy. 
This adds more confirmation to the observation of a singularity in the energy supercritical case with cubic nonlinearity $p=1$. 

In addition, we can do numerical 
simulations for $p=1$ and other supercritical choices of the fractional parameter $s<\frac{1}{4}$. Some results, using super-Gaussian initial data \eqref{ininonlin}(c), can be seen in the box below:
\begin{center}
\begin{tabular}{ |c| c| c| }
\hline
$s$ & $t_{\rm f}$ & $\mu$\\
\hline
0.2 & 5.46 & -0.32 \\ 
 0.19 & 5.014 &  -0.31\\  
 0.15 & 3.9354 & -0.29\\
  \hline    
\end{tabular}
\end{center}

The fact that in all of these cases the solution $u$ behaves qualitatively similar, suggests that the singularity does {\it not} depend on a particular choice 
of the parameter $s<s_\ast$, at least in the cubic case $p=1$. 
\begin{remark}
One might suspect from our findings that the parameter $\mu\in \R$ goes to zero as $s\to 0$. This would be consistent with the fact that the limiting equation \eqref{limeq}, formally obtained for $s=0$, does not suffer from a singularity. 
However, we were not able to give conclusive numerical evidence for a (possible) dependence of $\mu$ on $s\in (0,1]$. 
\end{remark}


\subsection{Other nonlinearities}
Finally, we shall try to answer the question whether the previous findings are sensitive to the choice of the nonlinear power $p\in \N$. 

To this end, we first consider a quintic nonlinearity with $p=2$. Here, the critical index is $s_\ast = \frac{1}{3}$, hence bigger than in the previous (cubic) case. More specifically, we shall study 
the supercritical time-evolution of
\begin{equation}\label{quinfNLS}
i \eps \partial_t u  = (- \eps^2 \Delta )^{1/4} u + |u|^{4} u, \quad u_{\vert{t=0}} = e^{-x^4}.
\end{equation}

\begin{figure}[htb!]
 \includegraphics[width=0.49\textwidth]{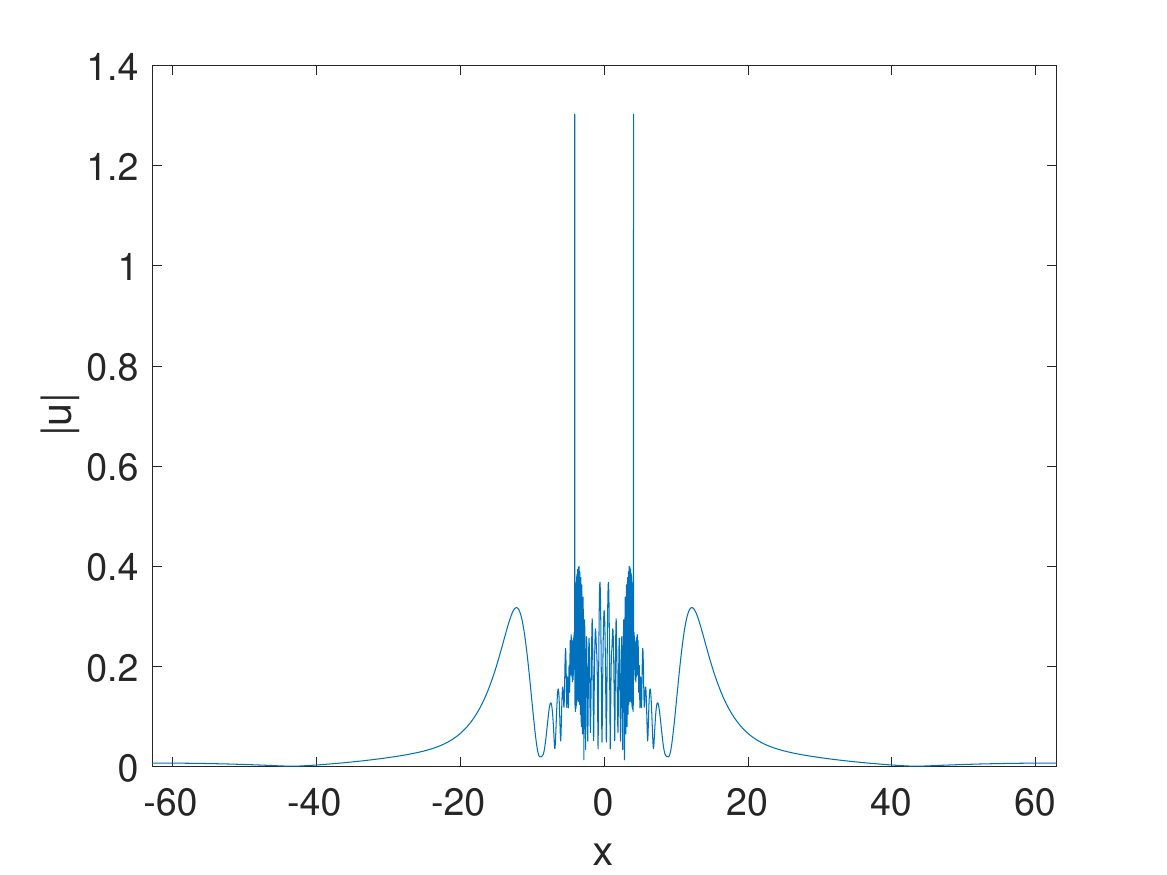}
 \includegraphics[width=0.49\textwidth]{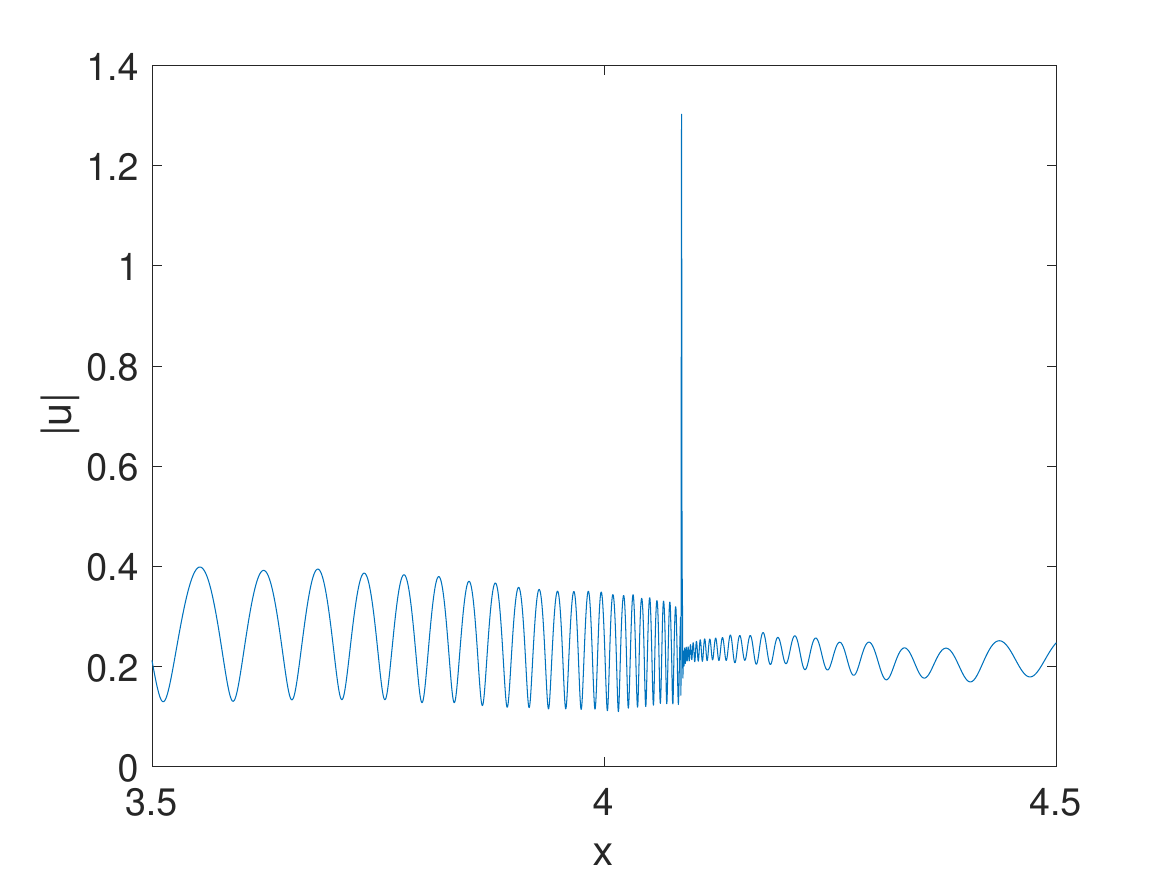}
 \caption{Left: solution to the quintic fNLS  
 \eqref{quinfNLS} with $\eps=0.1$ and $s=\frac{1}{4}$ at the final time $t_{\rm f}=8.47$. Right: a close-up 
 near one of the peaks.}
 \label{fNLSs02p2t847}
\end{figure}

Numerically we find $\delta (t) = 0$ at time $t_{\rm f} = 8.47$. The associated parameter $\mu = -0.38\dots$ and thus,
again compatible with $\mu \simeq -\frac13$. The solution 
at the final time is shown in Fig.~\ref{fNLSs02p2t847}. Overall, the quintic case behaves qualitatively similar to the one with a cubic nonlinearity. 
In particular, the oscillations again become strongly compressed near the presumed singularity at $x_0 \simeq \pm 4.1$.

\begin{figure}[htb!]
 \includegraphics[width=0.49\textwidth]{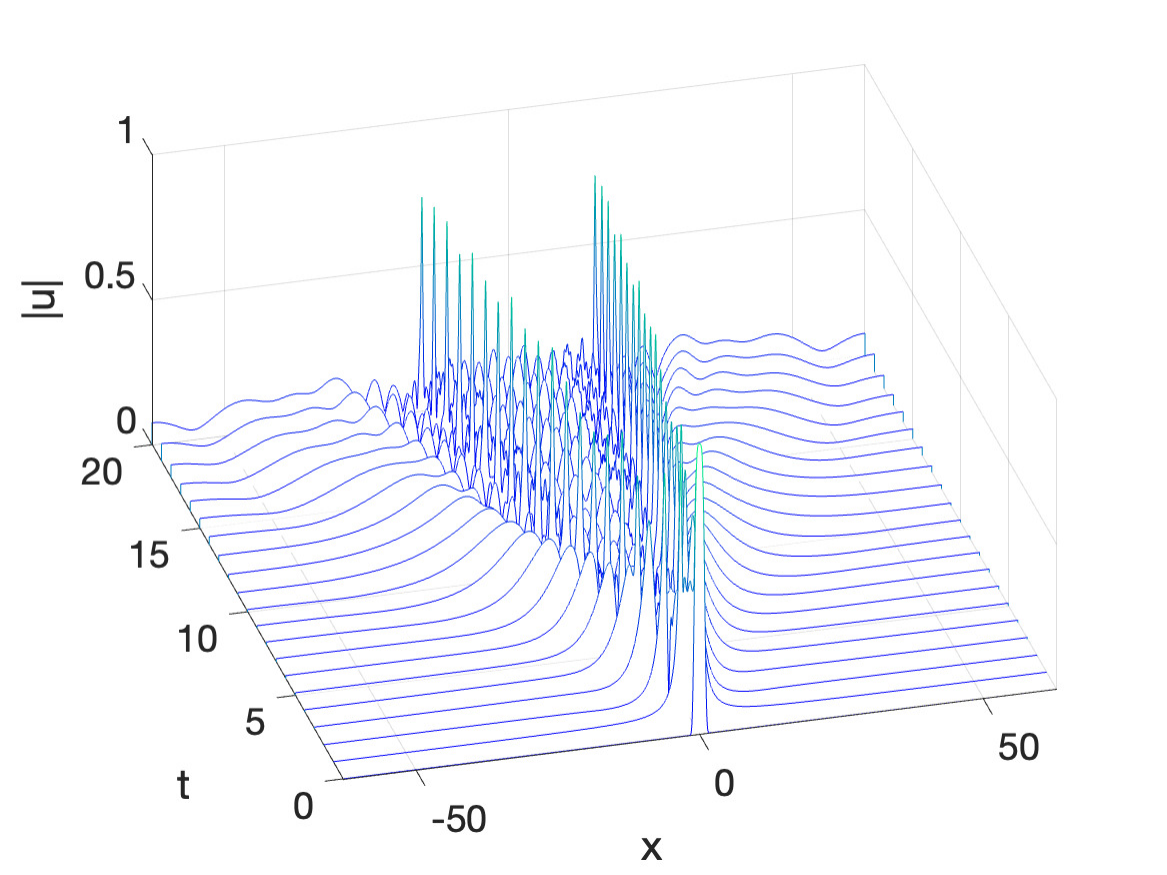}
 \includegraphics[width=0.49\textwidth]{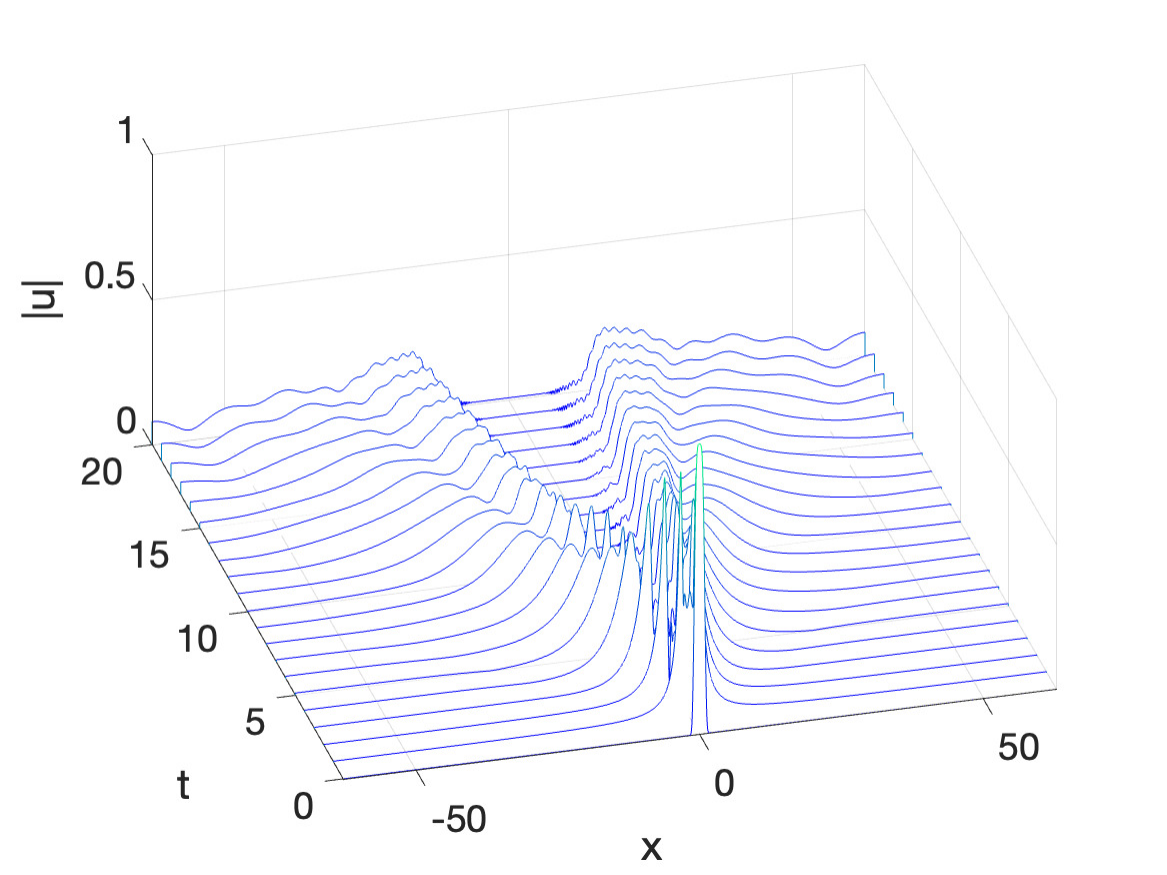}
 \caption{Solution to the fNLS equation \eqref{resfNLS} with $\eps=0.1$,
 $s=\frac{1}{4}$ and initial data \eqref{ininonlin}(c): on the left for 
nonlinear power $p=3$ and on the right for $p=4$.}
 \label{fNLSs0252p}
\end{figure}

Surprisingly, the same super-Gaussian initial data do not seem to yield a singularity for even higher nonlinear powers $p\ge 3$.  It can be seen 
on the left of Fig.~\ref{fNLSs0252p} that no singularity appears for $p=3$, even though the solution still displays a weakly turbulent behavior. 
Even more surprisingly, the right picture within Fig.~\ref{fNLSs0252p} shows that for $p=4$ the 
time-evolution yields a purely dispersive solution, similar to the linear case (compare with Fig. \eqref{fschr025}). 

\begin{figure}[htb!]
 \includegraphics[width=0.49\textwidth]{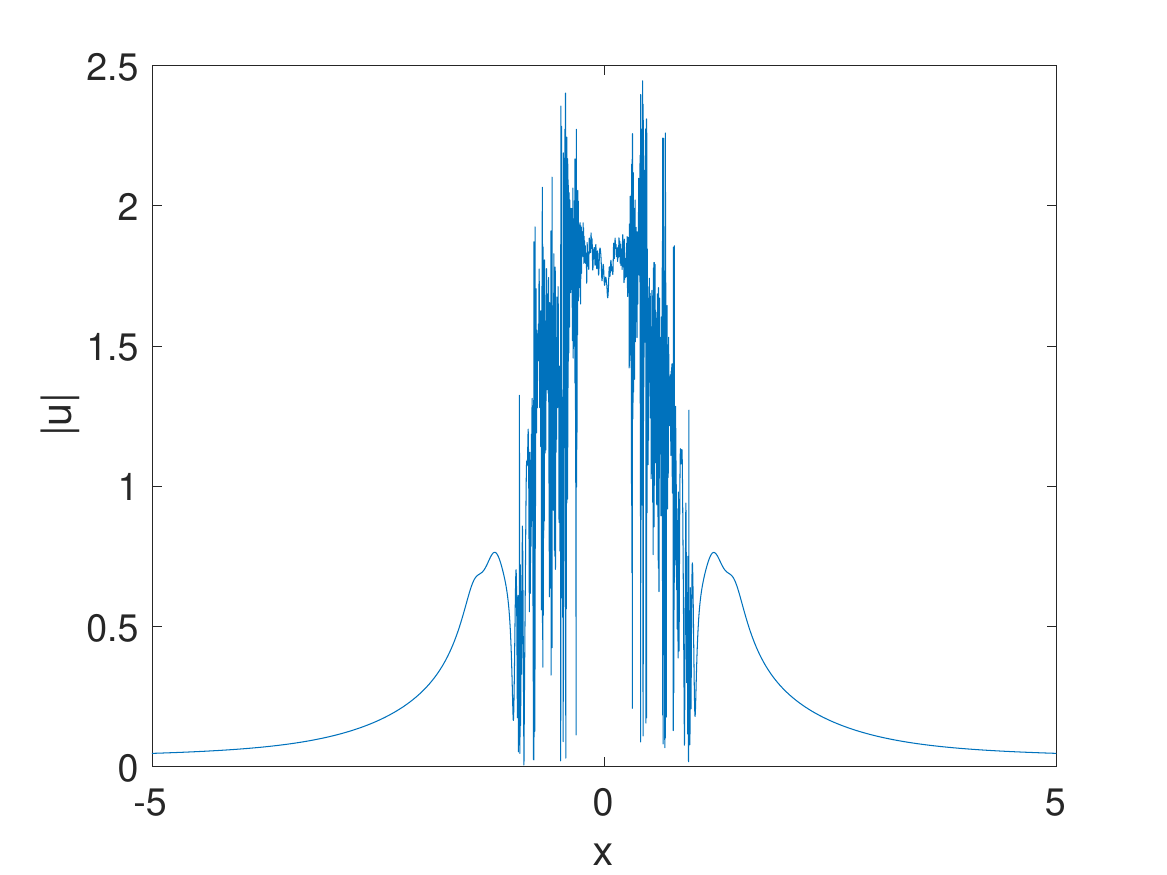}
 \includegraphics[width=0.49\textwidth]{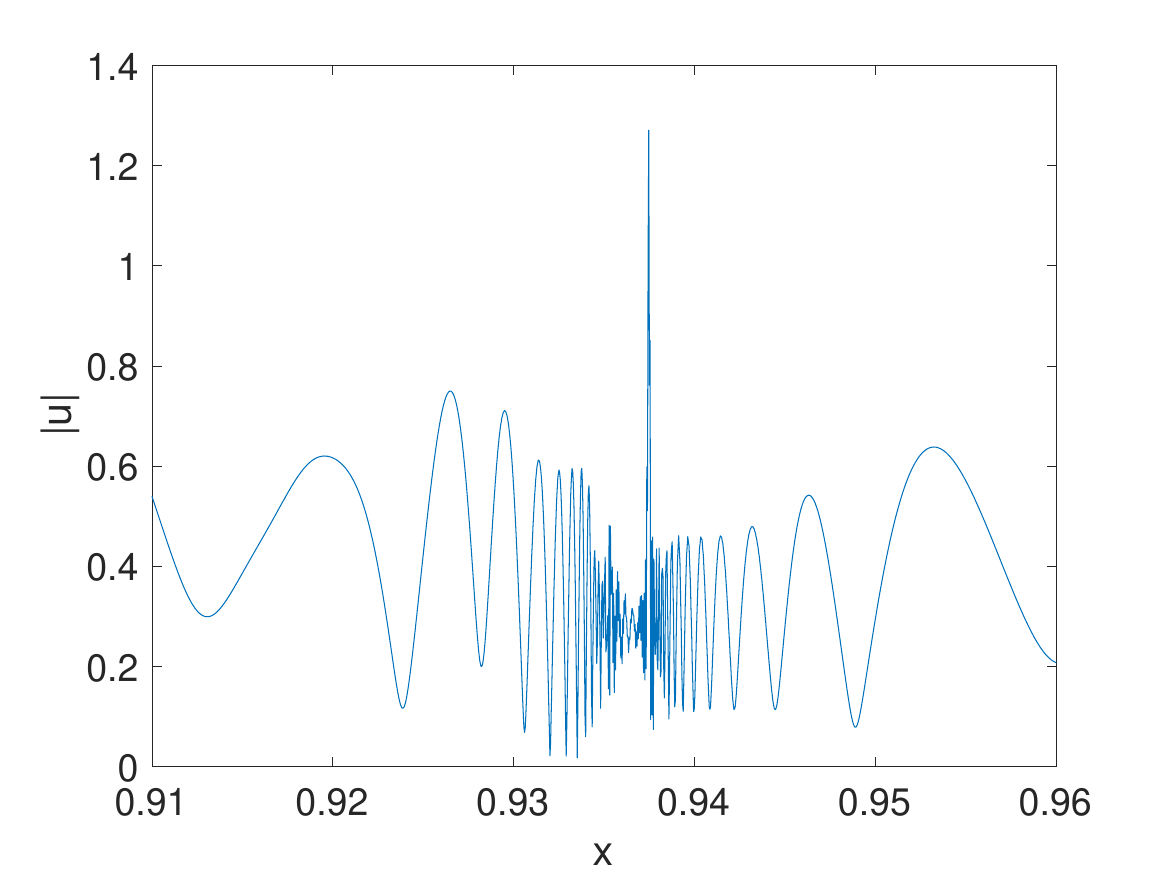}
 \caption{Left: Solution to the fNLS equation \eqref{resfNLS} with $\eps=0.1$,
 $s=\frac{1}{4}$, $p=3$, and initial data \eqref{2ini} at the final time $t_{\rm f}=1$. Right:
a close-up of the solution near one of the peaks.}
 \label{fNLSs025p3_2gauss}
\end{figure}

However, once we chose bigger initial data, such as 
\begin{equation}\label{2ini}
u_{\vert{t=0}}(x)=2e^{-x^{2}},
\end{equation}
a singular behavior again reemerges. Indeed, we find that in this case the code stops already at 
$t_{\rm f} = 1$, to yield a fitted $\mu\simeq -0.31$. The solution 
at the final time is shown on the left of 
Fig.~\ref{fNLSs025p3_2gauss}. A close-up near one of the peaks (shown at the right of the same 
figure) indicates a compression of the phase similarly to our previous supercritical cases. 

These findings suggest that such highly oscillatory singularities can only be observed for initial data with sufficiently large $L^q$-norms. 
It is beyond the scope of the current article to quantify the associated threshold and its precise dependence on the nonlinear power $p\in \N$. 
At this point, it is also unclear if the highly oscillatory dynamics produced in our numerical simulations is of the same type as the blow-up described in \cite{MRRS3}. 
Since our results strongly depend on the choice of $s<s_\ast$, it could be that our (presumed) singularities are of a different nature than the ones found in \cite{MRRS3}. 
As far as we know, however, there have not yet been any successful numerical simulations of the singularity formation proved in \cite{MRRS3}. 


\bigskip

\appendix

\section{Comparison with other numerical methods}

The convergence of the numerical scheme (RK4) used in the present work was discussed in detail in 
\cite{etna} for the case of the classical nonlinear Schr\"odinger equation ($s=1$). 
In the following we briefly compare and contrast (RK4) with some other 
numerical approaches one might want to pursue: 

We first note that the exponential integrators discussed in \cite{etna} 
can in principle also be applied in the present case, but with the following caveats: the 
exponential integrator appearing in our case is of the form
\[
S(h) = e^{-i\varepsilon^{2s-1}|k|^{2s}h},
\] 
where the time-step is $\Delta t = h>0$. Multiplying the numerical solution by such a factor becomes problematic 
whenever $\varepsilon^{2s-1}|k|^{2s}h=0\mbox{ mod }2\pi$, which 
implies the CFL-type condition $h<2\pi/(\varepsilon^{2s-1}|k|^{2s})$. 
On the one hand, this condition can enforce time steps that effectively become too small in our case which requires sufficiently high 
resolution (in Fourier space). On the other hand, Driscoll's \cite{Dri} composite Runge-Kutta method can suffer from 
instabilities for larger time steps which enforces less efficient 
implementations, see, e.g. the discussion in \cite{GKP}. 

\begin{figure}[htb!]
 \includegraphics[width=0.49\textwidth]{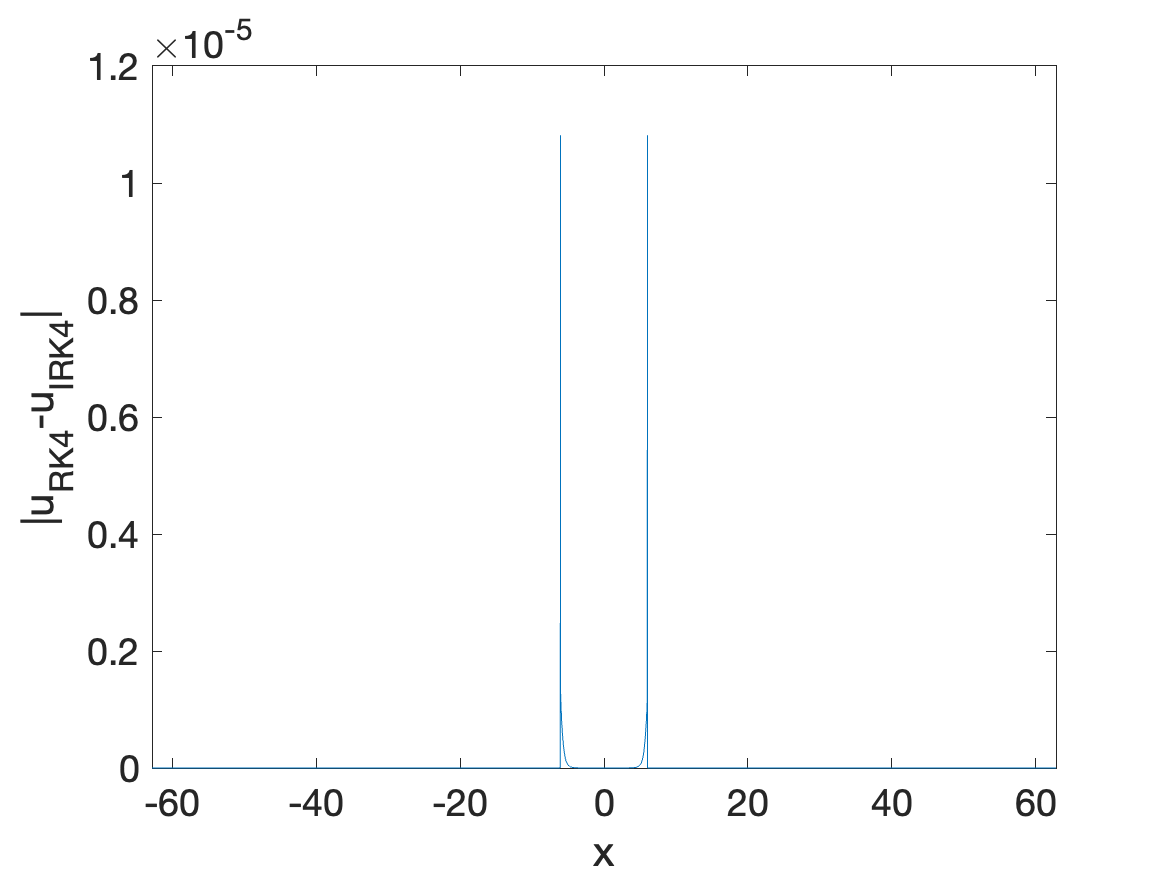}
 \includegraphics[width=0.49\textwidth]{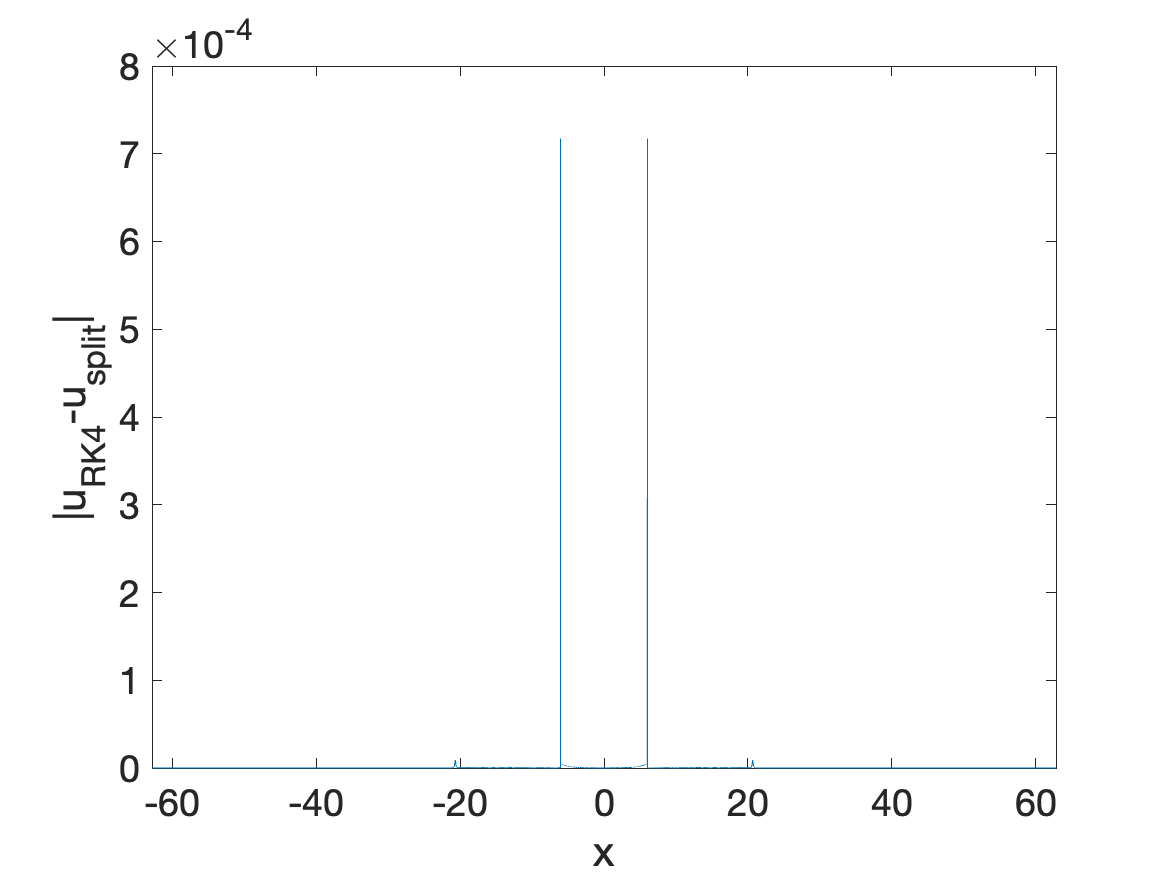}
 \caption{Left: Difference between the numerical solution of Figure 
 \ref{fNLSs025t20} and a solution obtained using IRK4. 
 Right: Difference between the numerical solution of Figure 
 \ref{fNLSs025t20} and a fourth order splitting method.}
 \label{fNLSs025t20deltat}
\end{figure}

The most efficient method from 
those discussed in \cite{etna} are fourth order splitting schemes, originally described in
\cite{Yos}. Below, we show a comparison of our (explicit) RK4 algorithm with such a fourth order splitting scheme, as well as with an implicit Runge-Kutta method (IRK4). 
We thereby focus on the example given in 
Fig.~\ref{fNLSs025t20}, i.e. the fNLS with cubic nonlinearity ($p=1$) in the energy-critical case ($s=\frac14$). We again chose $\eps = 0.1$ and initial data $\upsilon(x) = \text{sech}(x)$. 
In Fig.~\ref{fNLSs025t20deltat} we show the difference in $L^\infty$-norm between our results and the other two aforementioned methods. We thereby use the same numerical 
parameters $N$ and $N_t$ for all three methods, ensuring numerical errors smaller than $10^{-3}$. As expected, the resulting difference is of the order of 
$10^{-4}$, the biggest deviations appearing at the sharp peaks within the 
solution $u$. 

We note that the splitting scheme is roughly 30\% faster (Matlab 
timings) than our RK4. The reason why we only use it as reference 
solution in the present paper is that numerical resonances limit the 
achievable accuracy to the order $\mathcal O(10^{-8})$, cf. the discussion in \cite{etna}.

\bigskip

\bibliographystyle{amsplain}

\end{document}